\documentclass[11pt]{article}

\usepackage{amsmath}
\usepackage{amssymb}
\usepackage{latexsym}
\usepackage{comment}

\newtheorem{assumption}{Assumption}
\def\qed{ \ \vrule width.2cm height.2cm depth0cm\smallskip}
\newenvironment{proof}{\noindent {\bf Proof.\/}}{$\qed$\vskip 0.1in}

\newcommand{\ba}{\begin{array}}
\newcommand{\ea}{\end{array}}
\newcommand{\be}{\begin{equation}}
\newcommand{\ee}{\end{equation}}
\newcommand{\bea}{\begin{eqnarray}}
\newcommand{\eea}{\end{eqnarray}}
\newcommand{\beaa}{\begin{eqnarray*}}
\newcommand{\eeaa}{\end{eqnarray*}}

\def\dbE{\mathbb{E}}
\def\dbF{\mathbb{F}}

\def\dbL{\mathbb{L}}

\def\dbP{\mathbb{P}}
\def\dbR{\mathbb{R}}

\def\a{\alpha}
\def\b{\beta}
\def\g{\gamma}
\def\d{\delta}
\def\e{\varepsilon}

\def\t{\tau}
\def\f{\varphi}
\def\th{\theta}
\def\o{\omega}

\def\D{\Delta}

\def\L{\Lambda}

\def\O{\Omega}
\def\cC{{\cal C}}

\def\cE{{\cal E}}
\def\cF{{\cal F}}

\def\cP{{\cal P}}

\def\cT{{\cal T}}

\def\ch{\textsc{h}}

\def\no{\noindent}

\def\ms{\medskip}

\def\q{\quad}
\def\qq{\qquad}

\def\cd{\cdot}
\def\cds{\cdots}

\def\tr{\hbox{\rm tr}}

\def\qed{ \hfill \vrule width.25cm height.25cm depth0cm\smallskip}

\newcommand{\basa}{\begin{assumption}}
\newcommand{\easa}{\end{assumption}}

\newcommand{\bas}{\begin{assum}}
\newcommand{\eas}{\end{assum}}

\def\limsup{\mathop{\overline{\rm lim}}}
\def\liminf{\mathop{\underline{\rm lim}}}

 \def\cd{\cdot}
\def\cds{\cdots}

\def\tr{\hbox{\rm tr$\,$}}

\def\wh{\widehat}

\def\1{{\bf 1}}

\def\:{\!:\!}
\def\reff#1{{\rm(\ref{#1})}}
\def \proof{{\noindent \bf Proof\quad}}

\def \dbf{{\mathbf{d}}}

at 9pt

\begin{document}

\newtheorem{thm}{Theorem}[section]
\newtheorem{lem}[thm]{Lemma}
\newtheorem{cor}[thm]{Corollary}
\newtheorem{prop}[thm]{Proposition}
\newtheorem{rem}[thm]{Remark}
\newtheorem{eg}[thm]{Example}
\newtheorem{defn}[thm]{Definition}
\newtheorem{assum}[thm]{Assumption}

\renewcommand {\theequation}{\arabic{section}.\arabic{equation}}
\def\thesection{\arabic{section}}

\title{\bf Optimal Stopping under Nonlinear Expectation}

\author{Ibrahim  {\sc Ekren}\footnote{University of Southern California, Department of Mathematics, ekren@usc.edu.}   \and 
Nizar {\sc Touzi}\footnote{CMAP, Ecole Polytechnique Paris, nizar.touzi@polytechnique.edu. Research supported by the Chair {\it Financial Risks} of the {\it Risk Foundation} sponsored by Soci\'et\'e G\'en\'erale, and
              the Chair {\it Finance and Sustainable Development} sponsored by EDF and Calyon. }
       \and Jianfeng {\sc Zhang}\footnote{University of Southern California, Department of Mathematics, jianfenz@usc.edu. Research supported in part by NSF grant DMS 10-08873.}
}\maketitle

\begin{abstract}
Let $X$ be a bounded c\`adl\`ag process with positive jumps defined on the canonical space of continuous paths. We consider the problem of optimal stopping the process $X$ under a nonlinear expectation operator $\cE$ defined as the supremum of expectations over a weakly compact family of nondominated measures. We introduce the corresponding nonlinear Snell envelope. Our main objective is to extend the Snell envelope characterization to the present context. Namely, we prove that the nonlinear Snell envelope is an $\cE-$supermartingale, and an $\cE-$martingale up to its first hitting time of the obstacle $X$. This result is obtained under an additional uniform continuity property of $X$. We also extend the result in the context of a random horizon optimal stopping problem.

This result is crucial for the newly developed theory of viscosity solutions of path-dependent PDEs as introduced in \cite{EKTZ}, in the semilinear case, and extended to the fully nonlinear case in the accompanying papers \cite{ETZ1,ETZ2}.
\end{abstract}



\vfill\eject

\section{Introduction}
\label{sect-Introduction}
\setcounter{equation}{0}

On the canonical space of continuous paths, we consider a bounded c\`adl\`ag process $X$, with positive jumps, and satisfying some uniform continuity condition. Let $\ch_0$ be the first exit time of the canonical process from some convex domain, and $\ch:=\ch_0\wedge t_0$ for some $t_0>0$. This paper focuses on the problem
 \beaa
 \sup_{\t\in\cT}\cE[X_{\tau\wedge\ch}],
 &\mbox{where}&
 \cE[.]:=\sup_{\dbP\in\cP}\dbE^\dbP[.],
 \eeaa
$\cT$ is the collection of all stopping times, relative to the natural filtration of the canonical process, and $\cP$ is a weakly compact non-dominated family of singular measures.

Our main result is the following. Similar to the standard theory of optimal stopping, we introduce the corresponding nonlinear Snell envelope $Y$, and we show that the classical Snell envelope characterization holds true in the present context. More precisely, we prove that the Snell envelope $Y$ is an $\cE-$supermartingale, and an $\cE-$martingale up to its first hitting time $\tau^*$ of the obstacte. Consequently, $\tau^*$ is an optimal stopping time for our problem of optimal stopping under nonlinear expectation.

This result is proved by adapting the classical arguments available in the context of the standard optimal stopping problem under linear expectation. However, such an extension turns out to be highly technical. The first step is to derive the dynamic programming principle in the present context, implying the $\cE-$supermartingale property of the Snell envelope $Y$. To establish the $\cE-$martingale property on $[0,\tau^*]$, we need to use some limiting argument for a sequence $Y_{\t_n}$, where $\t_n$'s are stopping times increasing to $\t^*$. However, we face one major difficulty related to the fact that in a nonlinear expectation framework the dominated convergence theorem fails in general.  It was observed in Denis, Hu and Peng \cite{DHP} that the monotone convergence theorem holds in this framework if the decreasing sequence of  random variables are quasi-continuous. 
%convergence theorems under nonlinear expectation require some regularity of the sequence of random variables. This was observed in Denis, Hu and Peng \cite{DHP}. 
Therefore, one main contribution of this paper is to construct convenient quasi-continuous approximations of the sequence $Y_{\t_n}$. This allows us to apply the arguments in \cite{DHP} on $Y_{\t_n}$, which is decreasing under expectation (but not pointwise!) due to the supermartingale property. The weak compactness of the class $\cP$ is crucial for the limiting arguments.

We note that in an one dimensional Markov model with uniformly non-degenerate diffusion, Krylov \cite{Krylov} studied a similar optimal stopping problem in the language of stochastic control (instead of nonlinear expectation). However, his approach relies heavily on the smoothness of the (deterministic) value function, which we do not have here. Indeed,  one of the main technical difficulties in our situation is to obtain the locally  uniform regularity of the value process.

Our interest in this problem is motivated from the recent notion of viscosity solutions of path-dependent partial differential equations, as developed in \cite{EKTZ} and the accompanying papers \cite{ETZ1,ETZ2}. Our definition is in the spirit of Crandal, Ishii and Lions \cite{CIL}, see also Fleming and Soner \cite{FS}, but avoids the difficulties related to the fact that our canonical space fails to be locally compact. 
The key point is that the pointwise maximality condition, in the standard theory of viscosity solution, is replaced by a problem of optimal stopping under nonlinear expectation. 

Our previous paper \cite{EKTZ} was restricted to the context of semilinear path-dependent partial differential equations. In this special case, our definition of viscosity solutions can be restricted to the context where $\cP$ consists of absolutely continuous measures on the canonical space. Consequently, the Snell envelope characterization of the optimal stopping problem under nonlinear expectation is available in the existing literature on reflected backward stochastic differential equations, see e.g. El Karoui et al \cite{EKPPQ}, Bayraktar, Karatzas and Yao \cite{BayraktarKaratzas}. However, the extension of our definition to the fully nonlinear case requires to consider a nondominated family of singular measures. 

The paper is organized as follows. Section  \ref{sect-preliminary} introduces the probabilistic framework. Section \ref{sect-problemformulation} formulates the problem of optimal stopping under nonlinear expectation, and contains the statement of our main results. The proof of the Snell envelope characterization in the deterministic maturity case is reported in Section \ref{sect-T}. The more involved case of a random maturity is addressed in Section \ref{sect-ch}.

\section{Nondominated family of measures on the canonical space}
\label{sect-preliminary}
\setcounter{equation}{0}
 
\subsection{The canonical spaces}

Let $\O:= \big\{\o\in C([0,T], \dbR^d): \o_0={\bf 0}\big\}$, the set of continuous paths starting from the origin, $B$ the canonical process, $\dbF$ the natural filtration generated by $B$, 
$\dbP_0$ the Wiener measure,  $\cT$   the set of  $\dbF$-stopping times, and %$\L := \{(t, \o_{\cd\wedge t}): (t,\o)\in [0,T]\times \O\}$.
$\L := [0,T]\times \O$. Here and in the sequel, for notational simplicity,  we use ${\bf 0}$ to denote vectors or matrices with appropriate dimensions whose components are all equal to $0$. 
% For some technical convenience, we denote
 %\bea
 %\label{|x|}
 %|x| := \max_{1\le i\le d} |x_i| &\mbox{for}& x \in \dbR^d,
 %\eea
 %and then 
 %We define a norm on $\O$ and a metric on $\L$ as follows: 
 We define a seminorm on $\O$ and a pseudometric on $\L$ as follows: for any $(t, \o), ( t', \o') \in\L$, 
\bea\label{rho}
 \|\o\|_{t} 
 := 
 \sup_{0\le s\le t} |\o_s|,
 \q  
 \dbf_\infty\big((t, \o),( t', \o')\big) 
 := 
 |t-t'| + \big\|\o_{.\wedge t} - \o'_{.\wedge t'}\big\|_T.
 \eea
Then $(\O, \|\cd\|_{T})$ is a Banach space and $(\L, \dbf_\infty)$ is a complete pseudometric space. In fact, the subspace  $\{(t, \o_{\cd\wedge t}): (t,\o)\in \L\}$ is a complete metric space under $\dbf_\infty$.  
%For arbitrary  $(t,\o)\in [0,T]\times \O$,
%when there is no confusion, we may abuse the notation by identifying $(t, \o)$ with $(t, \o_{\cd\wedge t})$ and thus denote $(t,\o) \in \L$. 

\ms

We next introduce the shifted spaces. Let  $0\le s\le t\le T$.

-  Let $\O^t:= \big\{\o\in C([t,T], \dbR^d): \o_t ={\bf 0}\big\}$ be the shifted canonical space; $B^{t}$ the shifted canonical process on
$\O^t$;   $\dbF^{t}$ the shifted filtration generated by $B^{t}$, $\dbP^t_0$ the Wiener measure on $\O^t$,  $\cT^t$  the set of  $\dbF^t$-stopping times, and $\L^t :=  [t,T]\times \O^t$.

%- Define $\|\cd\|^t_s$, $\dbf^t_\infty$,  and $C^0(\L^t)$ etc. %, $C_b^0(\L^t)$, UC$(\L^t)$, and UC$_b(\L^t)$ 
%in the spirit of (\ref{rho}) and Definition \ref{defn-spaceC0}.

- For  $\o\in \O^s$ and $\o'\in \O^t$, define the concatenation path $\o\otimes_{t} \o'\in \O^s$ by:
\beaa
(\o\otimes_t \o') (r) := \o_r\1_{[s,t)}(r) + (\o_{t} + \o'_r)\1_{[t, T]}(r),
&\mbox{for all}&
r\in [s,T].
\eeaa

- Let $s\in[0,T)$ and $\o\in \O^s$. For an $\cF^{s}_{T}$-measurable random variable $\xi$, an $\dbF^{s}$-progressively measurable
process $X$ on $\O^s$, and $t\in(s,T]$, define the shifted $\cF^{t}_{T}$-measurable random variable $\xi^{t,\o}$  and  $\dbF^{t}$-progressively measurable
process $X^{t,\o}$ on $\O^t$ by:
\beaa
\xi^{t, \o}(\o') :=\xi(\o\otimes_t \o'), \q X^{t, \o}(\o') := X(\o\otimes_t \o'),
&\mbox{for all}&
\o'\in\O^t.
\eeaa

\subsection{Capacity and nonlinear expectation}

A probability measure $\dbP$ on $\O^t$ is called a semimartingale measure if the canonical process $B^t$ is a semimartingale under $\dbP$. For every constant $L>0$, we denote by $\cP^L_t$ the collection of all semimartingale measures $\dbP$ on $\O^t$ such that  there exist $\dbF^t$-progressively measurable $\dbR^d$-valued process $\a^\dbP$, a process $\b^\dbP\ge {\bf 0}$ with $d\times d$-symmetric matrix values, and a $d$-dimensional $\dbP$-Brownian motion $W^\dbP$ satisfying:
\bea
\label{cPL}
d B_t =\b^\dbP_tdW^\dbP_t + \a^\dbP_t dt,
&\dbP\mbox{-a.s. and}&
|\alpha^\dbP|\le L,~~\tr((\beta^\dbP)^2)\le 2L.
\eea
Throughout this paper, we shall consider a family $\{\cP_t,t\in[0,T]\}$ of semimartingale measures on $\O^t$ satisfying:
\begin{enumerate}
\item[(P1)] {\it there exists some $L_0$ such that, for all $t$,  $\cP_t$ is a weakly compact subset of $\cP_t^{L_0}$.}
\item[(P2)] {\it For any $0\le t\le T$, $\t\in\cT^t$, and $\dbP\in \cP_t$, the r.c.p.d. $\dbP^{\t,\o} \in \cP_{\t(\o)}$ for $\dbP$-a.e. $\o\in \O^t$.}
\item[(P3)] {\it For any  $0\le s\le t\le T$, $\dbP\in \cP_s$, $\{E_i, i\ge 1\} \subset \cF^s_t$ disjoint, and $\dbP^i\in \cP_t$, the following $\hat\dbP$ is also in $\cP_s$:}
\bea
\label{hatP}
\hat\dbP &:=& \dbP \otimes_t \Big[\sum_{i=1}^\infty \dbP^i \1_{E_i} + \dbP \1_{\cap_{i=1}^\infty E_i^c}\Big].
\eea
\end{enumerate}
Here \reff{hatP} means, for any event $E\in\cF^s_T$ and denoting $E^{t,\o} := \{\o'\in \O^t: \o\otimes_t \o' \in E\}$:
 \beaa
 \hat\dbP[E] 
 := 
 \dbE^\dbP\Big[\sum_{i=1}^\infty \dbP^i[E^{t,\o}] \1_{E_i}(\o) \Big]
 +\dbP\big[E\cap(\cap_{i=1}^\infty E_i^c)\big].
\eeaa
We refer to the seminal work of Stroock and Varadhan \cite{SV} for the introduction of r.c.p.d.,  which is a convenient tool for proving the dynamic programming principles, see e.g.  Peng \cite{Peng-g} and Soner, Touzi, and Zhang \cite{STZ-duality}.

We observe that for all $L>0$, the family $\{\cP_t^L,t\in[0,T]\}$ satisfies conditions (P1-P2-P3). In particular, the weak compactness follows standard arguments, see e.g. Zheng \cite{zheng} Theorem 3. The following are some other typical examples of such a family $\{\cP_t,t\in[0,T]\}$. 

\begin{eg}
\label{eg-cP} {\rm Let $L,  L_1, L_2> 0$ be some constants.

\no {\it Wiener measure}  $\cP^0_t := \{\dbP_0^t\} = \{\dbP: \a^\dbP=0, \b^\dbP = I_d\}$.

\no {\it Finite variation} $\cP^{\mbox{\sc fv}}_t(L) := \{\dbP:  |\a^\dbP|\le L, \b^\dbP={\bf 0}\}$.

\no {\it Drifted Wiener measure} $\cP^{0,{\rm ac}}_t(L) := \{\dbP:  |\a^\dbP|\le L,  \b^\dbP=I_d\}$.

\no {\it Relaxed bounds} $\cP_t(L_1,L_2) := \{\dbP: |\a^\dbP|\le L_1,  {\bf 0} \le \b^\dbP \le L_2 I_d\}$. 

\no {\it Relaxed bounds, Uniformly elliptic} $\cP_t^{\mbox{\sc ue}}(L_1,L_2,L)  :=  \{\dbP: |\a^\dbP|\le L_1,  L I_d\le   \b^\dbP \le L_2 I_d\}$.

\no {\it Equivalent martingale measures} 
$
\cP^{\rm e}_t(L_1,L_2,L)
:=  \{\dbP\in\cP_t(L_1,L_2)\!: \exists\,|\g^\dbP|\le L,
                               \a^\dbP \!=\! \b^\dbP \g^\dbP
    \}.
$
}
\end{eg}

We denote by $\dbL^1(\cF^t_T,\cP_t)$ the set of all $\cF^t_T-$measurable r.v. $\xi$ with $\sup_{\dbP\in\cP_t}\dbE^{\dbP}[|\xi|]<\infty$. The set $\cP_t$ induces the following capacity and nonlinear expectation:
  \bea\label{cE}
 \cC_t[A]
 :=
 \sup_{\dbP\in\cP_t}\dbP[A]
 ~
 \mbox{for}~
 A\in\cF^t_T, &\mbox{and}& 
 {\cE}_t[\xi]
 =
 \sup_{\dbP\in\cP_t}\dbE^{\dbP}[\xi]
 ~\mbox{for}~
 \xi\in\dbL^1(\cF^t_T,\cP_t).
 \eea
 When $t=0$, we shall omit $t$ and abbreviate them as $\cP, \cC, \cE$.
Clearly $\cE$ is a $G$-expectation, in the sense of Peng \cite{Peng-G}. We remark that, when $\xi$ satisfies certain regularity condition,  then $\cE_t[\xi^{t,\o}]$ can be viewed as the conditional $G$-expectation of $\xi$, and as a process it is the solution of a Second Order BSDEs, as introduced by Soner, Touzi and Zhang \cite{STZ-2BSDE}. 

Abusing the terminology of Denis and Martini \cite{DM}, we say that a property holds $\cP$-q.s. (quasi-surely) if it holds $\dbP-$a.s. for all $\dbP\in\cP$. A random variable $\xi: \O\to \dbR$ is 
\\
- $\cP$-quasicontinuous if for any $\e>0$,  there exists a closed set $\O_\e\subset \O$ such that $\cC(\O_\e^c) <\e$ and $\xi$ is continuous in $\O_\e$,
\\
- $\cP$-uniformly integrable if $\cE[|\xi|\1_{\{|\xi|\ge n\}}\longrightarrow 0$, as $n\to\infty$. 

Since $\cP$ is weakly compact, by Denis, Hu and Peng \cite{DHP} Lemma 4 and Theorems 22,28, we have:

\begin{prop}\label{prop-capacity}
{\rm (i)} Let $(\O_n)_{n\ge 1}$ be a sequence of open sets with $\O_n\uparrow\O$. Then $\cC(\O_n^c)\downarrow 0$.
\\
{\rm (ii)} Let $(\xi_n)_{n\ge 1}$ be a sequence of  $\cP$-quasicontinuous and $\cP$-uniformly integrable maps from $\O$ to $\dbR$. If $\xi_n\downarrow \xi$, $\cP$-q.s. then $\cE[\xi_n]\downarrow \cE[\xi]$. 
\end{prop}

We finally recall the notion of martingales under nonlinear expectation.

\begin{defn}
Let $X$ be an $\dbF$-progressively measurable process with $X_\t\in\dbL^1(\cF_\t,\cP)$ for all $\t\in \cT$. We say that $X$ is a ${\cE}-$supermartingale (resp. submartingale, martingale) if, for any $(t,\o)\in \L$ and any $\t\in\cT^t$, $\cE_t[X^{t,\o}_\t]\le$ (resp. $\ge,=$) $X_t(\o)$ for $\cP$-q.s. $\o\in\O$.
\end{defn}
We remark that we require the $\cE$-supermartingale property holds for stopping times. Under linear expectation $\dbP$, this is equivalent to the $\dbP$-supermartingale property for deterministic times, due to the Doob's optional sampling theorem. However, under nonlinear expectation, they are in general not equivalent.

\section{Optimal stopping under nonlinear expectations}
\label{sect-problemformulation}
\setcounter{equation}{0}

We now fix an $\dbF$-progressively measurable process $X$. 
 
\begin{assum}
\label{assum-X}
$X$ is a bounded c\`adl\`ag process with positive jumps, and there exists a modulus of continuity function $\rho_0$ such that for any $(t,\o), (t',\o')\in\L$:
 \bea\label{USC}
 X(t,\o) - X(t',\o') 
 \le 
 \rho_0\Big(\dbf_\infty\big((t,\o), (t',\o')\big)\Big)
 ~\mbox{whenever}~t\le t'.
 \eea 
\end{assum}

\begin{rem}\label{rem-USC}
{\rm There is some redundancy in the above assumption. Indeed, it is shown at the end of this section that {\it \reff{USC} implies that $X$ has left-limits and $X_{t-} \le X_t$ for all $t\in(0,T]$}. Moreover, the fact that $X$ has only positive jumps is important to ensure that the random times $\t^*$ in \reff{Yt*}, $\hat\t^*$ in \reff{hatY}, and $\t_n$ in  \reff{T-tn} and \reff{Hn} are $\dbF$-stopping times.\qed}
\end{rem}

We define the nonlinear Snell envelope and the corresponding obstacle first hitting time:
\bea
\label{Yt*}
Y_t (\o) := \sup_{\t\in\cT^t} \cE_t[X^{t,\o}_\t], 
&\mbox{and}& 
\t^* := \inf\{t\ge 0: Y_t = X_t\}.
\eea
Our first result is the following {\it nonlinear Snell envelope} characterization of the deterministic maturity optimal stopping problem $Y_0$.

\begin{thm}{{\rm (Deterministic maturity)}}
\label{thm-T-optimal}
Let $X$ be satisfying Assumption \ref{assum-X}.
Then $Y$ is an $\cE$-supermartingale on $[0,T]$,  $Y_{\t^*} = X_{\t^*}$, and $Y_{.\wedge\t^*}$ is an $\cE$-martingale. Consequently, $\t^*$ is an optimal stopping time for the problem $Y_0$. 
\end{thm}

To prove the partial comparison principle for viscosity solutions of path-dependent partial differential equations in our accompanying paper \cite{ETZ2}, we need to consider optimal stopping problems with random maturity time $\ch \in \cT$ of the form
 \bea\label{ch}
 \ch 
 := 
 \inf\{t\ge 0: B_t \in O^c\} \wedge t_0,
 \eea
for some $t_0\in(0,T]$ and some open convex set $O  \subset \dbR^d$ containing the origin. We shall extend the previous result to the following stopped process:
 \bea\label{hatX}
 \wh X^\ch_s := X_s \1_{\{s <\ch\}} + X_{\ch-}\1_{\{s\ge \ch\}}
 &\mbox{for}&
 s\in[0,T].
 \eea
The corresponding Snell envelope and obstacle first hitting time are denoted:
 \bea\label{hatY}
 \wh Y_t^\ch(\o):= \sup_{\t\in \cT^t}
 \cE_t\Big[\big(\wh X^\ch\big)^{t,\o}_\t\Big],
 &\mbox{and}&
 \wh\t^*:=\inf\{t\ge 0:\wh Y^\ch_t=\wh X^\ch_t\}.
 \eea
 
Our second main result requires the following additional assumption.

\begin{assum}\label{assum-additional}
{\rm (i)} For some $L>0$, $\cP^{\mbox{\sc fv}}_t(L)\subset \cP_t$ for all $t\in[0,T]$, where $\cP^{\mbox{\sc fv}} _t(L)$ is defined in Example \ref{eg-cP}. 

\no {\rm (ii)} For any $0\le t<t+\d\le T$, $\cP_t \subset \cP_{t+\d}$ in the following sense: for any $\dbP\in \cP_t$ we have $\tilde \dbP\in \cP_{t+\d}$, where $\tilde \dbP$ is the probability measure on $\O^{t+\d}$ such that the $\tilde\dbP$-distribution of $B^{t+\d}$ is equal to the $\dbP$-distribution of $\{B^t_{s}, t\le s\le T-\d\}$.
\end{assum} 

\begin{thm}{{\rm (Random maturity)}}\label{thm-optimal}
Let $X$ be a process satisfying Assumption \ref{assum-X}, and suppose that the nondominated family of singular measures satisfies Assumption \ref{assum-additional}. 
Then  $\wh Y^\ch$ is an $\cE$-supermartingale on $[0,\ch]$,  $\wh Y^\ch_{\wh\t^*} = \wh X^\ch_{\wh\t^*}$, and $\wh Y^\ch_{.\wedge\wh\t^*}$ is an $\cE$-martingale. In particular, $\wh\t^*$ is an optimal stopping time for the problem $\wh Y^\ch_0$. 
\end{thm}

\begin{rem}\label{rem-ch}
{\rm (i) The main idea for proving Theorem \ref{thm-optimal} is to show that $\cE[\wh Y^\ch_{\t_n}]$ converges to $\cE[\wh Y^\ch_{\wh \t^*}]$, where $\t_n$ is defined by \reff{Hn} below and increases to $\wh\t^*$. However, we face a major difficulty that the dominated convergence theorem fails in our nonlinear expectation framework. Notice that $Y$ is an $\cE$-supermartingale and thus $Y_{\t_n}$ are decreasing under expectation (but not pointwise!). We shall extend the arguments of \cite{DHP} for  the monotone convergence theorem, Proposition \ref{prop-capacity}, to our case. 
%The main difficulty for proving Theorem \ref{thm-optimal} is that  stopping times  are not continuous in $\o$, and thus we cannot apply the monotone convergence theorem, Proposition \ref{prop-capacity}, directly on them. 
For this purpose, we need to construct certain continuous approximations of the stopping times $\t_n$, and the requirement that the random maturity $\ch$ is of the form \reff{ch} is crucial.
% in order to construct a quasicontinuous approximation of $\wh\tau^*$.
We remark that, in his Markov model, Krylov \cite{Krylov} also considers this type of hitting times. We also remark that, in a special case,  Song \cite{Song} proved that $\ch$ is quasicontinuous.

\no (ii) Assumption \ref{assum-additional} is a technical condition used to prove the dynamic programming principle in Subsection 5.1 below. By a little more involved arguments, we may prove the results by replacing   Assumption \ref{assum-additional} (i)  with 
\beaa
\mbox{for some constant $L, L_1, L_2$},~~\cP_t^{\mbox{\sc ue}}(L_1,L_2,L) \subset \cP_t
&\mbox{for all}&
t\in[0,T],
\eeaa
where $\cP_t^{\mbox{\sc ue}} $ is defined in Example \ref{eg-cP} (iv).
\qed
}
\end{rem}

We conclude this section with the

\vspace{5mm}

\no {\bf Proof of Remark \ref{rem-USC}}\quad Fix $\o\in\O$, and let $ \{ t_n \} $ and $ \{s_n \} $ be two sequences such that $t_n \uparrow t, s_n \uparrow t$, and
$X_{t_n} \longrightarrow \limsup_{s \uparrow t} X_s$, $X_{s_n} \longrightarrow \liminf_{s \uparrow t} X_s$. Here and in the sequel, in $\lim_{s\uparrow t}$ we take the notational convention that $s<t$. 
Without loss of generality, we may assume $t_n < s_{n} <  t_{n+1}$ for $n = 1,2,...$.  Then for the $\rho_0$ defined in \reff{USC} we have 
$$
0\le  \limsup_{s \uparrow t} X_s
      -\liminf_{s \uparrow t} X_s
= \lim_{n \rightarrow \infty}X_{t_n}
  - \lim_{n \rightarrow \infty} X_{s_n}
\le
\liminf_{n\to\infty} 
\rho_0\Big(\dbf_\infty\big(t_n,\o), (s_n,\o)\big)\Big) 
=0.
$$
This implies the existence of $X_{t-}(\o)$. Moreover, 
$$
X_{t-} - X_{t} 
=  
\lim_{s \uparrow t} X_s -X_t 
\le  
\liminf_{s \uparrow t} \rho\Big(\dbf_\infty\big((s,\o), (t,\o)\big)\Big) =0,
$$
completing the proof.
\qed

\section{Deterministic maturity optimal stopping}
\label{sect-T}
\setcounter{equation}{0}

We now prove Theorem \ref{thm-T-optimal}. Throughout this section, Assumption \ref{assum-X} is always in force, and we consider the nonlinear Snell envelope $Y$ together with the first obstacle hitting time $\t^*$, as defined in \reff{Yt*}. Assume $|X|\le C_0$, and without loss of generality that $\rho_0\le 2C_0$. It is obvious that
\bea
\label{Ybasic}
\mbox{$|Y|\le C_0$, $Y\ge X$, and $Y_T = X_T$.}
\eea
Throughout this section, we shall use the following modulus of continuity function:
\bea
\label{rho0}
%\overline{|B^t|}_s:=\sup_{t\le r\le s}|B^t_r|,\q 
 \bar \rho_0(\d) := \rho_0(\d)\vee \Big[\rho_0(\d^{1\over 3}) + \d^{1\over 3}\Big],
\eea
and we shall use a generic constant $C$ which depends only on $C_0$, $T$, $d$, and the $L_0$ in Property (P1), and it may vary from line to line.

\subsection{Dynamic Programming Principle}

Similar to the standard Snell envelope characterization under linear expectation, our first step is to establish the dynamic programming principle. We start by the case of determinsitic times.

\begin{lem}
\label{lem-T-DPP1}
The process $Y$ is uniformly continuous in $\o$, with the modulus of continuity function $\rho_0$, and satisfies
 \bea\label{DP-det}
 Y_{t_1}(\o) 
 = 
 \sup_{\t\in \cT^{t_1}} 
 \cE_{t_1}
 \Big[X^{t_1,\o}_\t     \1_{\{\t  < t_2\}}
       + Y_{t_2}^{t_1,\o} 
         \1_{\{\t \ \ge t_2\}}
 \Big]
 &\mbox{for all}& 0\le t_1 \le t_2 \le T, \o\in \O.
 \eea
\end{lem}

\proof (i) First, for any $t$, any $\o, \o'\in \O$, and any $\t \in \cT^t$, by \reff{USC} we have
\beaa
|X^{t,\o}_\t - X^{t, \o'}_\t| &=& \Big|X(\t(B^t), \o\otimes_t B^t) - X(\t(B^t), \o'\otimes_t B^t)\Big| \\
&\le& \rho_0\Big(\dbf_\infty\big((\t(B^t), \o\otimes_t B^t), (\t(B^t), \o'\otimes_t B^t)\big)\Big) = \rho_0\big(\|\o-\o'\|_t\big).
\eeaa
Since $\t$ is arbitrary, this proves uniform continuity of $Y$ in $\o$.
\\
(ii) When $t_2=T$, since $Y_T=X_T$ \reff{DP-det} coincides with the definition of $Y$. Without loss of generality we assume $(t_1,\o)=(0,\mathbf{0})$ and $t:=t_2<T$. 
\\
{\it Step 1.} We first prove "$\le$". For any $\t\in \cT$ and $\dbP\in {\cP}$:
\beaa 
\dbE^ {\dbP }\left [X_{\t}\right]
&=&
\dbE^ {\dbP }\left [ X_{\t} \1_{\{\t < t\}} +\dbE_{t}^ {\dbP }[ X_{\t}] \1_{\{\t \ge t\}}\right]
\eeaa
By the definition of the r.c.p.d., we have $\dbE_{t}^ {\dbP }[X_{\t}] (\o)=\dbE^ {\dbP^{t,\o} }[X_{\t^{t,\o}}^{t,\o}]\leq Y_{t} (\o)$ for $\dbP-$a.e. $\o\in\{\t \ge t\}$, where the inequality follows from Property (P2) of the family $\{\cP_t\}$ that $\dbP^{t,\o} \in \cP_{t}$.
Then:
\beaa
\dbE^ {\dbP }\left [X_{\t}\right]
&\le&
\dbE^ {\dbP }\left [ X_{\t} \1_{\{\t < t\}} +Y_t \1_{\{\t \ge t\}}\right].
\eeaa
By taking the sup over $\t$ and $\dbP$, it follows that: 
 \beaa
 %\label{T-ineq1-DPPdeterministic}
 Y_0
 &=&
 \sup_{\t\in \cT} 
 \cE[ X_{\t }]
 \;\le\; 
 \sup_{\t\in \cT} 
 \cE\big[X_{\t }
                      \1_{\{\t < t\}}  
                      + Y_{t}
                        \1_{\{\t  \ge t\}} 
                \big].
 \eeaa
{\it Step 2.} We next prove "$\ge$". Fix arbitrary  $\t\in\cT$ and $\dbP\in\cP$, we shall prove
\bea
\label{T-DPP1est}
 \dbE^\dbP\big[X_{\t }
                      \1_{\{\t < t\}}  
                      + Y_{t}
                        \1_{\{\t  \ge t\}} 
                \big]&\le& Y_0.
\eea
Let $\e>0$, and  $\{E_i\}_{i\ge 1}$ be an $\cF_t$-measurable partition of  the event $\{\t\ge t\}\in\cF_t$ such that $\|\o-\tilde\o\|_t \le  \e$ for all $\o, \tilde\o \in E_i$. For each $i$,  fix an $\o^i\in E_i$, and  by the definition of $Y$ we have 
 \beaa
 Y_{t}(\o^{i})
 \;\le\;
 \dbE^{\dbP^{i}}
 \big[X^{t,\o^{i }}_{\t^{i}}
 \big]+\e
 &\mbox{for some}&
 (\t^{i},\dbP^{i})\in\cT^t\times\cP_t.
 \eeaa
By \reff{USC} and the uniform continuity of $Y$, proved in (i), we have
 \beaa
 |Y_t(\o) - Y_t(\o^i)|\le \rho_0(\e),\q |X^{t,\o}_{\t^{i}} - X^{t,\o^{i }}_{\t^{i}}| \le \rho_0(\e), &\mbox{for all}& \o\in E_i.
 \eeaa
Thus, for $\o\in E_i$,
  \bea
 \label{T-DPP1est1}
 Y_{t}(\o) \le Y_{t}(\o^{i}) + \rho_0(\e) \le 
 \dbE^{\dbP^{i}}
 \big[X^{t,\o^{i }}_{\t^{i}}
 \big]+\e + \rho_0(\e) \le  \dbE^{\dbP^{i}}
 \big[X^{t,\o}_{\t^{i}}
 \big]+\e + 2\rho_0(\e).
 \eea
Thanks to Property (P3) of the family $\{\cP_t\}$, we may define the following pair $(\tilde\t,\tilde\dbP)\in\cT\times\cP$:
\beaa
 \tilde\t :=\1_{\{\t< t\}} \t + \1_{\{\t\ge t\}}\sum_{i\ge 1}\1_{E_i} \t^{i}(B^t);&&
 \tilde \dbP :=\dbP \otimes_t \Big[ \sum_{i\ge 1} \1_{E_i}  \dbP^i +  \1_{\{\t < t\}}\dbP\Big].
 \eeaa
It is obvious that $\{\t<t\} = \{\tilde\t <t\}$. Then, by  \reff{T-DPP1est1},
 \beaa
 \dbE^{\dbP}
 \big[X_{\t}
      \1_{\{\t < t\}} 
      +Y_t \1_{\{\t \ge t\}}
 \big]
 &=& 
 \dbE^{\dbP}
 \Big[ X_{\t}
      \1_{\{\t < t\}} 
      +  \sum_{i\ge 1} Y_t \1_{E_i}
 \Big]\\
 &\le&   \dbE^{\dbP}
 \Big[ X_{\t}
      \1_{\{\t < t\}}   +  \sum_{i\ge 1}  \dbE^{\dbP^i} [ X^{t, \cd}_{\t^i}]  \1_{E_i} 
 \Big]+\e +2\rho_0(\e)\\
 &=& \dbE^{\tilde\dbP}
 \Big[X_{\tilde\t}
      \1_{\{\tilde\t < t\}} 
      + \sum_{i\ge 1} X_{\tilde\t} \1_{E_i}
 \Big]+\e +2\rho_0(\e)\\
 &=& \dbE^{\tilde\dbP}
 \big[ X_{\tilde\t} \big]+\e +2\rho_0(\e) \le Y_0 + \e +2\rho_0(\e),
 \eeaa
which provides \reff{T-DPP1est} by sending $\e\to 0$.  %Then the reverse inequality of \reff{T-ineq1-DPPdeterministic} follows from the  arbitrariness of $\dbP$ and $\t$.
\qed

\vspace{5mm}

We now derive the regularity of $Y$ in $t$. 

%\overline \rho_0(\d) := \sup_{0\le t\le T} \cE_t\Big[\rho_0\big(\d+ \sup_{t\le s\le (t+\d)\wedge T}|B^t_s|\big)\Big],
%\eea
%and observe that $\overline \rho_0(\d)\longrightarrow 0$ as $\d\searrow 0$, by the boundedness assumption in Property (P1) of the family $\cP_t$.

\begin{lem}
\label{lem-T-Ycontt}
For each $\o\in \O$ and  $0\le t_1<t_2\le T$, 
\beaa
|Y_{t_1}(\o) - Y_{t_2}(\o)| \le C\bar\rho_0\Big(\dbf_\infty\big((t_1,\o), (t_2,\o)\big)\Big).
\eeaa
\end{lem}
\proof Denote $\d := \dbf_\infty\big((t_1,\o), (t_2,\o)\big)$. If $\d \ge {1\over 8}$, then clearly $
|Y_{t_1}(\o) - Y_{t_2}(\o)|  \le 2C_0 \le C\bar \rho_0(\d)$. So we continue the proof assuming $\d \le {1\over 8} $.
First, by setting $\t = t_2$ in Lemma \ref{lem-T-DPP1},
\beaa
\d Y:=
Y_{t_2}(\o) -Y_{t_1}(\o) 
&\le& 
Y_{t_2}(\o)  - \cE_{t_1}\big[Y_{t_2}^{t_1,\o}\big] 
\\
&\le& 
\cE_{t_1}\big[ Y_{t_2}(\o)- Y_{t_2}(\o\otimes_{t_1} B^{t_1}) 
    \big] 
\\
&\le& 
 \cE_{t_1}\big[\rho_0\big(\dbf_\infty\big((t_2,\o), (t_2, \o\otimes_{t_1} B^{t_1})\big)\big)\big]
\\
&\le&
\cE_{t_1}\big[\rho_0\big(\d +  \|B^{t_1}\|_{t_1+\d}\big)\big]
.
\eeaa
On the other hand, by the inequality $X\le Y$, Lemma \ref{lem-T-DPP1}, and \reff{USC}, we have
\beaa
-\d Y
&\le&
\sup_{\t\in\cT^{t_1}}\cE_{t_1}
\Big[\big[X^{t_1,\o}_{t_2} 
          + 
          \rho_0\big(\dbf_\infty((\t,\o\otimes_{t_1}\! B^{t_1})
          ,(t_2, \o\otimes_{t_1}\!B^{t_1}))\big)
     \big] \1_{\{\t< t_2\}} \\
     &&\qq\qq\q
     + Y_{t_2}^{t_1,\o} \1_{\{\t\ge t_2\}}  \Big] -  Y_{t_2}(\o)\\
&\le& 
\cE_{t_1}\Big[ Y^{t_1,\o}_{t_2} - Y_{t_2}(\o) + \rho_0\big(\dbf_\infty((t_1, \o), (t_2, \o\otimes_{t_1} B^{t_1}))\big)\Big] \\
&\le& 
\cE_{t_1}\Big[ \rho_0\big(\dbf_\infty((t_2, \o), (t_2, \o\otimes_{t_1} B^{t_1}))\big) +\rho_0\big(\dbf_\infty((t_1, \o), (t_2, \o\otimes_{t_1} B^{t_1}))\big)\Big]
\\
&\le&  2\cE_{t_1}\big[ \rho_0\big(\d+\|B^{t_1}\|_{t_1+\d}\big)\big].
\eeaa
Hence 
$$
|\d Y|
\;\le\; 
2\cE_{t_1}\big[ \rho_0\big(\d+\|B^{t_1}\|_{t_1+\d}\big)\big]
\;\le\;
\cE_{t_1}\Big[\rho_0\big(\d +  {3\over 4} \d^{1\over 3}\big) + 2C_0 \1_{\{\|B^{t_1}\|_{t_1+\d}\ge  {3\over 4} \d^{1\over 3}\}}\Big].
$$ 
Since $\d + {3\over 4}\d^{1\over 3} \le \d^{1\over 3}$ for $\d \le {1\over 8}$, this provides:
\bea\label{rho1est}
|\d Y|
\;\le\;
\rho_0(\d^{1\over 3}) + C\d^{-{2\over 3}} \cE_{t_1}\Big[\|B^{t_1}\|_{t_1+\d}^2\Big] 
\;\le\; 
\rho_0(\d^{1\over 3}) + C\d^{-{2\over 3}}\d 
\le C\bar\rho_0(\d). 
\eea
\qed

We are now ready to prove the dynamic programming principle for stopping times. 
\begin{thm}\label{thm-T-optimal-DP}
For any $(t,\o)\in\L$ and $\t\in \cT^t$, we have 
 \beaa
 Y_t(\o) 
 =
 \sup_{\tilde\tau\in\cT^t}
\cE_t
 \Big[X^{t,\o}_{\tilde\t}
      \1_{\{\tilde\t< \t\}}
      +Y^{t,\o}_{\t}
       \1_{\{\tilde\t\ge \t\}}\Big].
 \eeaa
 Consequently, $Y$ is an $\cE$-supermartingale on $[0,T]$.
\end{thm}
\proof First, follow the arguments in Lemma \ref{lem-T-DPP1} (ii) Step 1 and note that Property (P2) of the family $\{\cP_t\}$ holds for stopping times, one can prove straightforwardly that
 \beaa
 Y_{t}(\o) 
 \le 
 \sup_{\tilde\t\in \cT^{t}} 
 \cE_t
 \Big[ X^{t,\o}_{\tilde\t}
       \1_{\{\tilde\t < \t\}}
       + Y_{\t}^{t,\o} 
         \1_{\{\tilde\t \ge \t\}}
 \Big].
\eeaa
On the other hand, let $\t_k \downarrow \t$ such that $\t_k$ takes only finitely many values. By Lemma \ref{lem-T-DPP1} one can easily show that Theorem \ref{thm-T-optimal-DP} holds for $\t_k$. Then for any $\dbP\in \cP_t$ and $\tilde\t\in \cT^t$, by denoting $\tilde\t_m := [\tilde\t + {1\over m}]\wedge T$ we have
 \beaa
 \dbE^\dbP\Big[X^{t,\o}_{\tilde\t_m}    \1_{\{\tilde\t_m< \t_k\}}
                + Y_{\t_k}^{t,\o} 
                  \1_{\{\tilde\t_m\ge\t_k\}}
          \Big] 
 \le Y_t(\o).
 \eeaa
Sending $k\to\infty$,  by Lemma \ref{lem-T-Ycontt} and the dominated convergence theorem (under $\dbP$):
 \beaa
 \dbE^\dbP\Big[X^{t,\o}_{\tilde\t_m}
               \1_{\{\tilde\t_m  \le \t\}}
               + Y_{\t}^{t,\o} 
                 \1_{\{\tilde\t_m>\t\}}
          \Big] 
 &\le& Y_t(\o).
 \eeaa
Since the process $X$ is right continuous in $t$, we obtain by sending $m\to \infty$:
 \beaa
 Y_t(\o)
 \ge
 \dbE^\dbP\Big[ X^{t,\o}_{\tilde\t }
               \1_{\{\tilde\t <\t\}}
               + Y_{\t}^{t,\o} 
                 \1_{\{\tilde\t\ge\t\}}
          \Big],
% \ge  \dbE^\dbP\Big[ X^{t,\o}_{\tilde\t }
 %              \1_{\{\tilde\t \le \t\}}
  %             + Y_{\t}^{t,\o} 
   %              \1_{\{\tilde\t>\t\}}
   %       \Big],
 \eeaa
%where the last inequality follows from the fact that $Y\ge X$. 
which provides the required result by the arbitrariness of $\dbP$ and $\tilde\t$.
\qed

\subsection{Preparation for the $\cE-$martingale property}

If $Y_0=X_0$, then $\t^*=0$ and obviously all the statements of Theorem \ref{thm-T-optimal} hold true. Therefore, we focus on the non-trivial case $Y_0>X_0$.

We continue following the proof of the Snell envelope characterization in the standard linear expectation context. Let
\bea
\label{T-tn}
\t_n := \inf\{t\ge 0: Y_t - X_t \le {1\over n}\} \wedge T
,
&\mbox{for}& n > (Y_0-X_0)^{-1}.
\eea

\begin{lem}
\label{lem-T-DPn}
The process $Y$ is an $\cE$-martingale on $[0,\t_n]$.
 \end{lem}
 \proof By the dynamic programming principle of Theorem \ref{thm-T-optimal-DP},
 \beaa
 Y_0 
 &=& 
 \sup_{\t\in \cT} 
 \cE\Big[X_{\t}\1_{\{\t < \t_n\}} 
                      +Y_{\t_n}\1_{\{\t\ge \t_n\}}
                 \Big].
 \eeaa
For any $\e>0$, there exist $\t_\e\in \cT$ and $\dbP_\e\in\cP$ such that
 \bea  
 \label{T-tauest1}
 Y_0 
\; \le\; 
 \dbE^{\dbP_\e}\Big[X_{\t_\e}
                      \1_{\{\t_\e< \t_n\}} 
                      + Y_{\t_n}\1_{\{\t_\e\ge \t_n\}}
                 \Big] 
 +\e
 &\le&
\dbE^{\dbP_\e}\Big[ Y_{\t_\e\wedge \t_n} 
                       -{1\over n}
                        \1_{\{\t_\e < \t_n\}}
                 \Big]
 +\e,
 \label{T-tauest2}
 \eea
where we used the fact that $Y_t - X_t > {1\over n}$ for $t<\t_n$, by the definition of $\tau_n$. 
On the other hand,  it follows from the $\cE-$supermartingale property of $Y$ in Theorem \ref{thm-T-optimal-DP} that $\dbE^{\dbP_\e}\Big[ Y_{\t_\e\wedge \t_n} \Big] \le \cE[Y_{\t_\e\wedge \t_n}]  \le Y_0$, which implies by \reff{T-tauest2} that $\dbP_\e[\t_\e < \t_n]\le n\e$. We then get from \reff{T-tauest1} that:
 $$
 Y_0 
 \le\dbE^{\dbP_\e}\!\Big[( X_{\t_\e} - Y_{\t_n})
                      \1_{\{\t_\e<\t_n\}} 
                      + Y_{\t_n}
                    \Big]  +\e
 \le 
 C\dbP_\e[\t_\e < \t_n] 
 + \dbE^{\dbP_\e}[Y_{\t_n}] +\e \le  \cE[Y_{\t_n}] + (Cn+1)\e
 . 
 $$
Since $\e$ is arbitrary, we obtain $Y_0 \le \cE[Y_{\t_n}]$. Similarly one can prove $Y$ is an $\cE$-submartingale on $[0, \t_n]$. By the $\cE-$supermartingale property of $Y$ established in Theorem \ref{thm-T-optimal-DP}, this implies that $Y$ is an $\cE-$martingale on $[0,\t_n]$.
\qed

\vspace{5mm}

By Lemma \ref{lem-T-Ycontt} we have
\bea
\label{T-Yt*}
Y_0 - \cE[Y_{\t^*}]  = \cE[Y_{\t_n}]- \cE[Y_{\t^*}] \le C\cE\Big[\bar\rho_0\Big(\dbf_\infty\big((\t_n, \o), (\t^*,\o)\big)\Big)\Big].
\eea
Clearly, $\t_n \nearrow \t^*$, and  $\bar\rho_0\Big(\dbf_\infty\big((\t_n, \o), (\t^*,\o)\big)\Big)\searrow 0$. However, in general the stopping times $\t_n, \t^*$ are not $\cP$-quasicontinuous, so we cannot apply Proposition \ref{prop-capacity} (ii) to conclude $Y_0 \le \cE[Y_{\t^*}]$. 
To overcome this difficulty, we need to approximate $\t_n$ by continuous r.v. 

\subsection{Continuous approximation}
 
The following lemma can be viewed as a Lusin theorem under nonlinear expectation and is crucial for us. 
\begin{lem}
\label{lem-qscont}
Let $\underline\th \le \th \le \overline\th$ be r.v. on $\O$, with values in a compact interval $I\subset\dbR$, such that for some $\O_0\subset \O$ and $\d>0$: 
 \beaa
 \underline\th(\o)  \le \th(\o') \le \overline \th(\o)   
 &\mbox{for all}&  
 \o\in \O_0 
 ~~\mbox{and}~~
 \|\o - \o'\| \le \d.
 \eeaa
Then for any $\e>0$, there exists a uniformly continuous function $\hat\th: \O \to I$ and an open subset $\O_\e \subset \O$ such that 
 \beaa
 \cC\big[\O_\e^c\big]\le\;\e\;
 &\mbox{and}& 
 \underline \th - \e \le \hat\th\le \overline \th +\e
 ~~\mbox{in}~~ \O_\e\cap \O_0.
 \eeaa
 \end{lem}

\proof If $I$ is a single point set, then $\th$ is a constant and the result is obviously true. Thus at below we assume the length $|I|>0$. Let $\{\o_j\}_{j\ge 1}$ be a dense sequence in $\O$. Denote $O_j := \{\o\in\O: \|\o-\o_j\|<{\d\over 2}\}$ and $\O_n := \cup_{j=1}^n O_j$. It is clear that $\O_n$ is open and $\O_n \uparrow \O$ as $n\to\infty$. Let $f_n : [0,\infty)\to [0,1]$ be defined as follows: $f_n(x) = 1$ for $x\in [0, {\d\over 2}]$, $f_n(x) = {1\over  n^2 |I|}$ for $x\ge \d$, and $f_n$ is linear in $[{\d\over 2}, \d]$. Define
 $$
 \th_n (\o) 
 := \phi_n(\o)\sum_{j=1}^n \th(\o_j) \f_{n,j}(\o)
 ~\mbox{where}~
 \f_{n,j}(\o):=f_n(\|\o-\o_j\|)
 ~\mbox{and}~
 \phi_n:= \Big(\sum_{j=1}^n  \f_{n,j}\Big)^{-1}.
 $$ 
Then clearly $\th_n$ is uniformly continuous and takes values in $I$. For each $\o \in \O_n\cap \O_0$, the set $J_n(\o) := \{ 1\le j\le n:  \|\o-\o_j\|\le \d\}\neq\emptyset$ and $\phi_n(\o)\le 1$. Then, by our assumption,
\beaa
 \th_n(\o)  - \overline \th (\o) 
 &=&  
 \phi_n(\o)\Big(\sum_{j\in J_n(\o)}
                       [\th(\o_j)-\overline\th(\o)] 
                       \f_{n,j}(\o) 
                + \sum_{j\notin J_n(\o)} 
                       [ \th(\o_j) - \overline \th(\o)]    
                       \f_{n,j}(\o)
           \Big)
 \\
 &\le&
 \phi_n(\o)\sum_{j\notin J_n(\o)}|I|\f_{n,j}(\o) 
 \;\le\; 
 \phi_n(\o)\sum_{j\notin J_n(\o)}{1\over n^2}
 \;\le\; {1\over n}.
 \eeaa
Similarly one can show that $\underline\th-{1\over n} \le \th_n$ in $\O_n\cap\O_0$. Finally, since  $\O_n \uparrow \O$ as $n\to\infty$, it follows from Proposition \ref{prop-capacity} (i)  that $\lim_{n\to\infty} \cC[\O_n^c] =0$. 
\qed

\subsection{Proof of Theorem \ref{thm-T-optimal}}  

We proceed in two steps.

\no {\it Step 1.} For each $n$, let $\d_n>0$ be such that $3C\bar\rho_0(\d_n) \le {1\over n(n+1)}$ for the constant $C$ in Lemma \ref{lem-T-Ycontt}. Now for any $\o$ and $\o'$ such that $\|\o-\o'\|_T \le \d_n$, by \reff{USC}, the uniform continuity of $Y$ in Lemma \ref{lem-T-DPP1}, and the fact that $\rho_0\le\bar\rho_0$, we have
\beaa
(Y-X)_{\t_{n+1}(\o)}(\o') 
\le 
(Y-X)_{\t_{n+1}(\o)}(\o) 
+ 3C\bar\rho_0(\d_n) 
\le 
{1\over n+1} + {1\over n(n+1)} 
= {1\over n}.
\eeaa 
Then $\t_n(\o') \le \t_{n+1}(\o)$. Since  $3C\bar\rho_0(\d_n) \le {1\over n(n+1)}\le {1\over n(n-1)}$, similarly we have $\t_{n-1}(\o) \le \t_{n}(\o')$. We may then apply Lemma \ref{lem-qscont} with $\underline{\theta}=\tau_{n-1}$, $\theta=\tau_n$,  $\overline{\theta}=\tau_{n+1}$, and $\O_0 = \O$. Thus, there exist an open set $\O_n\subset \O$ and a continuous r.v. $\tilde\t_n$ valued in $[0,T]$ such that
\beaa
\cC\big[\O_n^c\big] \le 2^{-n} &\mbox{and}& \t_{n-1} - 2^{-n} \le \tilde \t_n \le \t_{n+1} + 2^{-n}~\mbox{in}~\O_n.
\eeaa

\no {\it Step 2.}  By Lemma \ref{lem-T-DPn}, for each $n$ large, there exists $\dbP_n\in \cP$ such that 
\beaa
Y_0 = \cE[Y_{\t_n}] \le \dbE^{\dbP_n}[Y_{\t_n}] + 2^{-n}.
\eeaa
By Property (P1), $\cP$ is weakly compact. Then, there exists a subsequence $\{n_j\}$ and $\dbP^*\in \cP$ such that $\dbP_{n_j}$ converges weakly to $\dbP^*$. Now for any $n$ large and any $n_j\ge n$, note that $\t_{n_j} \ge \t_n$. Since $Y$ is an $\cE$-supermartingale and thus a $\dbP_{n_j}$-supermartingale,   we have 
 \bea
 \label{T-Y0e}
 Y_0 - 2^{-n_j} 
 &\le& 
 \dbE^{\dbP_{n_j}}\big[Y_{\t_{n_j}}\big] 
 \;\le\; 
 \dbE^{\dbP_{n_j}}\big[Y_{\t_{n}}\big] \le  \dbE^{\dbP_{n_j}}\big[Y_{\tilde \t_{n}}\big]  + \dbE^{\dbP_{n_j}}\big[|Y_{\tilde \t_{n}} - Y_{\t_n}|\big] .
 \eea
 By the boundedness of $Y$ in \reff{Ybasic} and the uniform continuity of $Y$ in Lemma \ref{lem-T-Ycontt}, we have
 \beaa
 |Y_{\tilde \t_{n}} - Y_{\t_n}| 
 &\le& 
 C\bar\rho_0\Big(\dbf_\infty\big((\tilde\t_n, \o), (\t_n, \o)\big)\Big)
 \\
 &\le& 
 C\bar\rho_0\Big(\dbf_\infty\big((\tilde\t_n, \o), (\t_n, \o)\big) \Big)\1_{\O_{n-1}\cap \O_{n+1}} + C \1_{\O_{n-1}^c \cup \O_{n+1}^c}.
 \eeaa
 Notice that $\tilde \t_{n-1} - 2^{1-n} \le \t_n \le \tilde\t_{n+1} + 2^{-1-n}$ on  $\O_{n-1}\cap \O_{n+1}$. Then
 \beaa
 |Y_{\tilde \t_{n}} - Y_{\t_n}| 
 &\le& 
 C\bar\rho_0
 \Big(\dbf_\infty\big((\tilde\t_n, \o), 
                   (\tilde \t_{n-1} - 2^{1-n}, \o)
              \big) \Big) \1_{\O_{n-1}\cap \O_{n+1}}
 \\
 &&
 +C\bar\rho_0
  \Big(\dbf_\infty\big((\tilde\t_n, \o), 
                    (\tilde \t_{n+1} +2^{-1-n}, \o)
               \big) \Big) \1_{\O_{n-1}\cap \O_{n+1}}
 + C \1_{\O_{n-1}^c \cup \O_{n+1}^c}
 \\
 &\le& 
 C\bar\rho_0
 \Big(\dbf_\infty\big((\tilde\t_n, \o), 
                   (\tilde \t_{n-1} - 2^{1-n}, \o)
              \big)\Big)  
 \\
 &&
 +C\bar\rho_0
 \Big(\dbf_\infty\big((\tilde\t_n, \o), 
                   (\tilde\t_{n+1}+2^{-1-n},\o)
              \big)\Big)  
 + C \1_{\O_{n-1}^c \cup \O_{n+1}^c}.
 \eeaa
Then \reff{T-Y0e} together with the estimate $\cC[\O_n^c]\le 2^{-n}$ lead to
 \beaa
 Y_0 - 2^{-n_j} 
 &\le&  
 \dbE^{\dbP_{n_j}}\big[Y_{\tilde \t_{n}}\big]  
 +C\dbE^{\dbP_{n_j}}
    \Big[\bar\rho_0
         \Big(\dbf_\infty\big((\tilde\t_n,\o), 
                           (\tilde \t_{n-1}-2^{1-n},\o)
                      \big)\Big)
    \Big] \\
 && 
 +C \dbE^{\dbP_{n_j}}
    \Big[\bar\rho_0
         \Big(\dbf_\infty\big((\tilde\t_n,\o), 
                           (\tilde \t_{n+1}+2^{-1-n},\o)
                      \big)\Big)
    \Big] 
 + C2^{-n}.
\eeaa
Notice that $Y$ and $\tilde\t_{n-1}, \tilde\t_n, \tilde\t_{n+1}$ are continuous. Send $j\to\infty$, we obtain
 \bea
 \label{T-Y0e2}
 Y_0 
 &\le&  
 \dbE^{\dbP^*}\big[Y_{\tilde \t_{n}}\big]  
 +C \dbE^{\dbP^*}
    \Big[\bar\rho_0
         \Big(\dbf_\infty\big((\tilde\t_n, \o),
                           (\tilde \t_{n-1}-2^{1-n},\o)
                      \big)\Big)
    \Big]
 \nonumber \\
 && 
 +C \dbE^{\dbP^*}
    \Big[\bar\rho_0
         \Big(\dbf_\infty\big((\tilde\t_n,\o),
                           (\tilde \t_{n+1}-2^{-1-n},\o)
                      \big)\Big)
    \Big] 
 + C2^{-n}.
\eea
Since $\sum_n \dbP^*\big[|\tilde\t_n - \t_n|\ge 2^{-n}\big] \le \sum_n \cC\big[|\tilde\t_n - \t_n|\ge 2^{-n}\big] \le \sum_n 2^{-n} <\infty$ and $\t_n \uparrow \t^*$, by the Borel-Cantelli lemma under $\dbP^*$
% in the capacity framework, see \cite{DHP} Lemma 9, and conclude that $\tilde\t_n \to \t^*$, $\cP$-q.s.  In particular, 
we see that $\tilde\t_n \to \t^*$, $\dbP^*$-a.s. Send $n\to\infty$ in \reff{T-Y0e2} and apply the dominated convergence theorem under $\dbP^*$, we obtain
\beaa
 Y_0 \le  \dbE^{\dbP^*}\big[Y_{\t^*}\big]  \le \cE[Y_{\t^*}].
 \eeaa
Similarly $Y_t(\o) \le   \cE_t[Y_{\t^*}^{t,\o}]$ for $t< \t^*(\o)$. By the $\cE$-supermartingale property of $Y$ established in Theorem \ref{thm-T-optimal-DP}, this implies that $Y$ is an $\cE$-martingale on $[0, \t^*]$.
\qed

\section{Random maturity optimal stopping}
\label{sect-ch}
\setcounter{equation}{0}

In this section, we prove Theorem \ref{thm-optimal}.  The main idea follows that of Theorem \ref{thm-T-optimal}. However, since $\wh X_\ch$ is not continuous in $\o$, the estimates become much more involved.

Throughout this section, let $X$, $\ch$, $O$, $t_0$, $\wh X:=\wh X^\ch$, $\wh Y:=\wh Y^\ch$, and $\wh\t^*$ be as in Theorem \ref{thm-optimal}. Assumptions  \ref{assum-X} and  \ref{assum-additional} will always be in force. We shall  emphasize when the additional Assumption \ref{assum-additional} is needed, and we fix the constant $L$ as in Assumption \ref{assum-additional}  (i).  
%However, unless stated explicitly, we do not assume $\wh \t^* = \ch$ in the discussion. 
Assume $|X|\le C_0$, and without loss of generality that $\rho_0 \le 2C_0$ and $L\le 1$. It is clear that 
\bea
\label{whYbasic}
\mbox{$|\wh Y|\le C_0$, $\wh X \le \wh Y$, and $\wh Y_{\ch} = \wh X_{\ch} = X_{\ch-}$.}
\eea
%Let $M_0$  denote the bound  of $X$.  Then $|\wh Y| \le M_0$, and 
%Without loss of generality we assume that the nondecreasing function $\rho_0$ is also bounded. 
By \reff{USC} and the fact that $X$ has positive jumps, one can check straightforwardly that,
 \bea 
 \label{USC1}
 &\wh X(t,\o) - \wh X(t', \o') \le  \rho_0\big(\dbf_\infty((t,\o), (t', \o'))\big)
 ~\mbox{for}~ t\le t',~t\le\ch(\o),~t'\le\ch(\o')&
 \\
% &\mbox{except the case}~t' = \ch(\o') <t_0~\mbox{and}~ t<\ch(\o).&\nonumber
 &\mbox{except the case}~ t=t' = \ch(\o') <\ch(\o)\le t_0.&\nonumber
 \eea 
%recall that \reff{USC} implies that $X_{t_-}\le X_t$ by Remark \ref{rem-USC}. 
In particular,
 \bea 
 \label{USC2}
 \wh X(t,\o) - \wh X(t', \o) \le  \rho_0\big(\dbf_\infty((t,\o), (t', \o))\big) &\mbox{whenever}& t\le t' \le \ch(\o).
 \eea 
 Moreover, we define
  \bea
\label{rho12}
%\overline{|B^t|}_s:=\sup_{t\le r\le s}|B^t_r|,\q 
 \rho_1(\d) :=  \rho_0(\d) \vee \Big[\rho_0\big((L^{-1}\d)^{1\over 3}\big)+ \d^{1\over 3}\Big] ,\q  \rho_2(\d) := [\rho_1(\d) + \d] \vee [\rho_1(\d^{1\over 3}) + \d^{1\over 3}],
\eea
 and in this section, the generic constant $C$ may depend on $L$ as well.
  %:=\max_{t\le r\le s}|B^t_r|$, and %$\Delta^{\!\!\eps}B_{[0,T]}:=\sup\{|B_t-B_s|:0\le s\le t\le T,t-s\le\eps\}$,
 %\bea\label{rhon}
 %C\rho_1(\d) 
 %&:=&  \sup_{0\le t\le T}
 %\cE_t 
 %\Big[\rho_0\Big(2L \d
  %                      + \sup_{t\le s\le T} |B^t_s- B^t_{(s-\d)\vee t}|\Big)\Big],
 %\\
 %\rho_n(\d) 
 %&:=&
 %\cE 
 %\Big[\rho_0\big(2\d + 2\overline{|B|}_\d\big)\Big] 
 %+ \cE
  % \Big[C\rho_1\big(\d + \overline{|B|}_\d\big)\Big]
  % + 2M_0\;\cC\Big[\overline{|B|}_\d \ge {1\over n} \Big],
  % \nonumber\\
 %\wh\rho_n(\d) 
 %&:=&
 %(2M_0 n \d) \vee C\rho_1(\d)\vee \rho_n(\d).\nonumber
 %\eea
%As $\d\searrow 0$,  by Property (P1) of $\cP_t$, it is clear that $\bar\rho_0(\d)\searrow 0$ and, $\rho_n(\d)\searrow 0$, $\bar\rho_n(\d)\searrow 0$ .

\subsection{Dynamic programming principle} 
%\label{sect-DPP}
%\setcounter{equation}{0}

We start with the regularity in $\o$.
\begin{lem}\label{lem-YUComega}
For any $t < \ch(\o) \wedge \ch(\o')$ we have: 
 \beaa
 |\wh Y_t(\o) - \wh Y_t(\o')| \le C\rho_1\big(\|\o-\o'\|_t\big).
\eeaa
\end{lem}
To motivate our proof, we first follow the arguments in Lemma \ref{lem-T-DPP1} (i) and see why it does not work here. Indeed, note that
\beaa
\wh Y_t(\o) - \wh Y_t(\o') \le \sup_{\t\in \cT^t} \sup_{\dbP \in \cP_t} \dbE^\dbP\Big[\wh X_{\t \wedge \ch^{t,\o}}^{t,\o} - \wh X_{\t \wedge \ch^{t,\o'}}^{t,\o'}\Big].
\eeaa 
Since we do not have  $\ch^{t,\o}\le \ch^{t,\o'}$, we cannot apply \reff{USC1} to obtain the required estimate.

\ms

\proof Let $\t \in \cT^t$ and $\dbP\in\cP_t$. Denote $\d := {1\over L}\|\o-\o'\|_t$,
% If  $t_0 - t \le \d$, then
%\beaa
 %|\wh Y_t(\o) - \wh Y_t(\o')| \le 2C_0 \le {C\d\over t_0-t} \le {C\rho_1(L\d)\over t_0-t}.
 %\eeaa
%From now on we assume $t_0-t >\d$.  Denote
 $t_\d := [t+\d]\wedge t_0$ and $\tilde B^{t_\d}_s := B^t_{s+\d} - B^t_{t_\d}$ for $s\ge t$.   Set   $\t'(B^{t}) := [\t(\tilde B^{t_\d}) +\d]\wedge t_0$, then $\t' \in \cT^t$. Moreover, by Assumption \ref{assum-additional} and Property (P3), we may choose $\dbP'\in  \cP_t$ defined as follows:  $\a^{\dbP'} := {1\over \d} (\o_{t}-\o'_{t})$, $\b^{\dbP'} :={\bf 0}$ on $[t, t_\d]$, and the $\dbP'$-distribution of $\tilde B^{t_\d}$  is equal to the $\dbP$-distribution of $B^t$. We claim that
\bea
\label{Ycontest1}
I:= \dbE^{\dbP} [\wh X^{t, \o}_{\t\wedge \ch^{t,\o}}]  - \dbE^{\dbP'} [\wh X^{t, \o'}_{\t'\wedge \ch^{t,\o'}}]  
&\le&  C\rho_1(L\d),
\eea
Then $
 \dbE^{\dbP} [\wh X^{t, \o}_{\t\wedge \ch^{t,\o}}] 
 - \wh Y_{t}(\o')  
 \le
 \dbE^{\dbP} [\wh X^{t, \o}_{\t\wedge \ch^{t,\o}}] 
 - \dbE^{\dbP'} [\wh X^{t, \o'}_{\t'\wedge \ch^{t,\o'}}]  \;\le\; C\rho_1(L\d)$, and it follows from the arbitrariness of $\dbP\in {\cP}_t$ and $\t \in \cT^t$ that $\wh Y_t(\o) - \wh Y_{t}(\o') \le   C\rho_1(L\d)$. By exchanging the roles of $\o$ and $\o'$, we obtain the required estimate.

It remains to prove \reff{Ycontest1}.  Denote
\beaa
\tilde\o'_s &:=& \o'_s\1_{[0,t)}(s) + [\o'_t + \a^{\dbP'}(s-t)]\1_{[t,T]}(s).
 %(\tilde\o'\otimes_{t_\d} B^{t})_s &:=& \tilde\o'_s\1_{[0,t_\d)}(s) + [\tilde\o'_{t_\d}+ B^t_{s-\d}]\1_{[t_\d, t_0]}(s).
\eeaa
%We note that $\tilde\o'\otimes_{t_\d} B^{t}$ involves a time shift of $B^t$ and is slightly different from the standard concatenation. 
Since $t<\ch(\o)\wedge\ch(\o')$, we have $\o_t, \o'_t \in O$. By the convexity of $O$, this implies that $\tilde\o'_s \in O ~\mbox{for}~s\in [t, t_\d]$, and thus 
%$ \ch^{t,\o'} > t_\d$ and 
$\ch^{t,\o'} (B^t)
 = (\ch^{t,\o}(\tilde B^{t_\d}) + \d)\wedge t_0,$ $\dbP'-$a.s. Therefore, 
 \bea
 \label{Ycontest2}
 \dbE^{\dbP'} [\wh X^{t, \o'}_{\t'\wedge \ch^{t,\o'}}] 
 &=& 
 \dbE^{\dbP'} \Big[\wh X\Big(\t'(B^{t})\wedge\ch^{t,\o'}(B^{t}),
                           \o'\otimes_{t} B^{t}\Big)\Big] \\
 &=& 
 \dbE^{\dbP'} \Big[\wh X\Big([\t(\tilde B^{t_\d})+\d]\wedge [\ch^{t,\o}(\tilde B^{t_\d}) +\d]\wedge t_0, \tilde\o'\otimes_{t_\d} \tilde B^{t_\d}_{\cd-\d}\Big)\Big]
 \nonumber\\
 &=&
 \dbE^{\dbP} \Big[\wh X\Big([\t(B^{t}) +\d]\wedge  
 [\ch^{t,\o}(B^{t}) +\d]\wedge t_0,  
 \tilde\o'\otimes_{t_\d} B^{t}_{\cd-\d}\Big)\Big],\nonumber
 \eea
while
 \beaa
 \dbE^{\dbP} [\wh X^{t, \o}_{\t\wedge \ch^{t,\o}}] 
 &=& 
 \dbE^{\dbP} \Big[\wh X\Big(\t(B^{t})\wedge\ch^{t,\o}(B^{t}),
                            \o\otimes_{t} B^{t}
                       \Big)
             \Big]. 
 \eeaa
 Notice that, whenever  $\t(B^{t})\wedge\ch^{t,\o}(B^{t}) = [\t(B^{t}) +\d]\wedge  
 [\ch^{t,\o}(B^{t}) +\d]\wedge t_0$, we have  $\t(B^{t})\wedge\ch^{t,\o}(B^{t})  = t_0$. This 
%In order to apply \reff{USC1}, we need to exclude the exceptional case there. In the case 
%\beaa
%t_0 &>& [\t(B^{t}) +\d]\wedge [\ch^{t,\o}(B^{t}) +\d]\wedge t_0 \\
%&=& \ch^{t, \tilde \o'}((\tilde\o'\otimes_{t_\d} B^{t}_{\cd-\d})^t) = \ch^{t, \o'}((\tilde\o'\otimes_{t_\d} B^{t}_{\cd-\d})^t) = [\ch^{t,\o}(B^t) +\d] \wedge t_0,
%\eeaa
%we have
%\beaa
%\ch^{t,\o}(B^t) +\d \le \t(B^{t}) +\d < t_0 &\mbox{and thus}& \ch^{t,\o}(B^t) \le \t(B^{t})  < t_0.
%\eeaa
%So $\t(B^t) \wedge \ch^{t,\o}(B^t) = \ch^{t,\o}(B^t)$, which 
excludes the exceptional case in \reff{USC1}. Then it  follows from \reff{Ycontest2} and \reff{USC1} that
\beaa
I
&\le&  \dbE^{\dbP} \Big[\rho_0\Big(\d +  \|(\o\otimes_{t} B^{t})_{\cd\wedge \t(B^t) \wedge \ch^{t,\o}(B^t)} - (\tilde\o'\otimes_{t_\d} B^{t}_{\cd-\d})_{\cd\wedge[\t(B^{t}) +\d]\wedge [\ch^{t,\o}(B^{t}) +\d]\wedge t_0} \|_{t_0} \Big)\Big].
\eeaa
Note that, denoting $\th := \t(B^t) \wedge \ch^{t,\o}(B^t)$,
 \beaa
&& \|(\o\otimes_{t} B^{t})_{\cd\wedge \t(B^t) \wedge \ch^{t,\o}(B^t)} - (\tilde\o'\otimes_{t_\d} B^{t}_{\cd-\d})_{\cd\wedge[\t(B^{t}) +\d]\wedge [\ch^{t,\o}(B^{t}) +\d]\wedge t_0} \|_{t_0}  \\
&\le& \|\o\otimes_{t} B^{t}- \tilde\o'\otimes_{t_\d} B^{t}_{\cd-\d} \|_{t_0} + \sup_{0\le r\le \d}|(\o\otimes_t B^t)_{ \th +r}- (\o\otimes_{t} B^{t})_{ \th}| \\
& \le&
\Big[ \|\o-\o'\|_t\Big] \vee \Big[ \sup_{t\le s\le t_\d}|\o_t + B^t_s - \tilde\o'_s|\Big] 
 \vee \Big[\sup_{t_\d \le s\le t_0}|\o_t + B^t_s 
                 - \tilde \o'_{t_\d} - B^t_{s-\d}|\Big]\\
                 &&+ \sup_{0\le r\le \d}|(\o\otimes_t B^t)_{ \th +r}- (\o\otimes_{t} B^{t})_{ \th}| \\
& \le& 2L\d +\|B^t\|_{t_\d} + \sup_{t_\d\le s\le t_0} |B^t_s - B^t_{s-\d}| + \sup_{0\le r\le \d}|B^t_{ \th +r}- B^{t}_{ \th}| .
 \eeaa
 Since $L\le 1$, we have
 \beaa
  I  &\le&   \dbE^{\dbP} \Big[\rho_0\Big(3\d + \|B^t\|_{t_\d} + \sup_{t_\d\le s\le t_0} |B^t_s - B^t_{s-\d}| + \sup_{0\le r\le \d}|B^t_{ \th +r}- B^{t}_{ \th}|\Big)\Big].
\eeaa
If $\d \ge {1\over 8}$, then $I \le 2C_0 \le C\rho_1(L\d)$. We then continue assuming $\d\le {1\over 8}$, and thus $3\d + {1\over 4} \d^{1\over 3} \le \d^{1\over 3}$. Therefore,
\beaa
I&\le& \rho_0(\d^{1\over 3}) + C\dbP\Big( \|B^t\|_{t_\d} + \sup_{t_\d\le s\le t_0} |B^t_s - B^t_{s-\d}| + \sup_{0\le r\le \d}|B^t_{ \th +r}- B^{t}_{ \th}| \ge {1\over 4} \d^{1\over 3}\Big)\\
&\le& \rho_0(\d^{1\over 3}) + C\d^{-{8\over 3}} \dbE^{\dbP}\Big[\|B^t\|^8_{t_\d} + \sup_{t_\d\le s\le t_0} |B^t_s - B^t_{s-\d}|^8 + \sup_{0\le r\le \d}|B^t_{ \th +r}- B^{t}_{ \th}|^8\Big]\\
&\le& \rho_0(\d^{1\over 3}) + C\d^{4\over 3} + C\d^{-{8\over 3}}\dbE^{\dbP}\Big[ \sup_{t_\d\le s\le t_0} |B^t_s - B^t_{s-\d}|^8\Big].
\eeaa
Set $t_\d = s_0 <\cds<s_n=t_0$ such that $\d \le s_{i+1}-s_i \le 2\d$, $i=0,\cds, n-1$. Then
\beaa
&&\dbE^{\dbP}\Big[ \sup_{t_\d\le s\le t_0} |B^t_s - B^t_{s-\d}|^8\Big] = \dbE^{\dbP}\Big[ \max_{0\le i\le n-1} \sup_{s_i \le s\le s_{i+1}}   |B^t_s - B^t_{s-\d}|^8\Big]\\
&\le&  \sum_{i=0}^{n-1}\dbE^{\dbP}\Big[ \sup_{s_i \le s\le s_{i+1}}   [|B^t_s - B^t_{s_i-\d}|+|B^t_{s-\d} - B^t_{s_i-\d}|]^8\Big]\\
&\le& C \sum_{i=0}^{n-1} (s_{i+1}-s_i + \d)^4 \le C\d^{-1} \d^4 = C\d^3.
\eeaa
Thus $I \le \rho_0(\d^{1\over 3}) + C\d^{4\over 3} + C\d^{-{8\over 3}} \d^3 \le  \rho_0(\d^{1\over 3}) + C\d^{1\over 3}   \le C\rho_1(L\d)$, proving \reff{Ycontest1} and hence the lemma.
\qed

\vspace{5mm}
We next show that the dynamic programming principle holds along deterministic times. 
\begin{lem}\label{lem-DPP1}
%Let Assumption \ref{assum-additional} hold true with constant $L$. Then, 
Let $t_1 < \ch(\o)$ and $t_2\in [t_1, t_0] $. We have: 
 \beaa
 \wh Y_{t_1}(\o) 
 &=& 
 \sup_{\t\in \cT^{t_1}} 
 \cE_{t_1}
 \Big[\wh X^{t_1,\o}_{\t \wedge\textsc{h}^{t_1,\o}}
       \1_{\{\t \wedge \textsc{h}^{t_1,\o} < t_2\}}
       + \wh Y_{t_2}^{t_1,\o} 
         \1_{\{\t \wedge \textsc{h}^{t_1,\o} \ge t_2\}}
 \Big].
 \eeaa
\end{lem}

\proof  When $t_2= t_0$, the lemma coincides  with the definition of $\wh Y$. Without loss of generality we assume $(t_1,\o)=(0,\mathbf{0})$ and $t := t_2<t_0$. First, follow the arguments in Lemma \ref{lem-T-DPP1} (ii) Step 1, one can easily prove  
 \bea\label{ineq1-DPPdeterministic}
 \wh Y_0
 &\le&
 \sup_{\t\in \cT} 
 \cE \Big[\wh X_{\t \wedge\textsc{h}}
                      \1_{\{\t \wedge \textsc{h} < t\}}  
                      + \wh Y_{t}
                        \1_{\{\t \wedge \textsc{h} \ge t\}} 
                \Big].
 \eea
To show that equality holds in the above inequality, fix arbitrary $\dbP\in\cP$ and $\t\in\cT$ satisfying $\t \le \ch$ (otherwise reset $\t$ as $\t\wedge \ch$),  we shall prove
\beaa
%\label{DPP1est}
\dbE^{\dbP}
 \Big[\wh X_{\t}
      \1_{\{\t < t\}} 
      +\wh Y_t \1_{\{\t \ge t\}}
 \Big] &\le& \wh Y_0.
\eeaa
Since $\wh Y_\ch = \wh X_\ch$, this amounts to show that:
\bea
\label{DPP1est}
\dbE^{\dbP}
 \Big[\wh X_{\t}
      \1_{\{\t< t\}\cup\{\ch \le t\}} 
      +\wh Y_t \1_{\{\t \ge t, \ch >t\}}
 \Big] &\le& \wh Y_0.
\eea
We adapt the arguments in Lemma \ref{lem-T-DPP1} (ii) Step 2 to the present situation.
%Let $L$ be such that Assumption \ref{assum-additional} holds. $\cP\subset\cP^L$ as in Property (P1), 
Fix   $0< \d \le  t_0 - t$. Let $\{E_i\}_{i\ge 1}$ be an $\cF_t$ measurable partition of the event $\{\t\ge t, \ch > t\} \in \cF_t$ such that $\|\o-\tilde\o\| \le L \d$ for all $\o, \tilde\o \in E_i$. Fix an $\o^i\in E_i$ for each $i$. By the definition of $\wh Y$ we have 
 \bea
 \label{DPP1est1}
 \wh Y_{t}(\o^{i})
 \;\le\;
 \dbE^{\dbP^{i}}
 \Big[\wh X^{t,\o^{i }}_{\t^{i} \wedge\ch^{t,\o^{i}}}
 \Big]+\d
 &\mbox{for some}&
 (\t^{i},\dbP^{i})\in\cT^t\times\cP_t.
 \eea
As in Lemma \ref{lem-YUComega}, we set $t_\d := t+\d<t_0$, $\tilde B^{t_\d}_s := B^t_{s+\d} - B^t_{t_\d}$ for $s\ge t$, and $\tilde \t^i(B^{t}) := [\t^i(\tilde B^{t_\d}) +\d]\wedge t_0$. Then $\tilde\t^i \in \cT^t$.  Moreover by Assumption \ref{assum-additional} and Property (P3), for each $\o\in E_i$, we may define $\dbP^{i,\o}\in  \cP_t$ as follows:  $\a^{\dbP^{i,\o}} := {1\over \d} (\o^i_{t}-\o_{t})$, $\b^{\dbP^{i,\o}} :={\bf 0}$ on $[t, t_\d]$, and the $\dbP^{i,\o}$-distribution of $\tilde B^{t_\d}$  is equal to the $\dbP^i$-distribution of $B^t$. By \reff{Ycontest1}, we have
\bea
\label{DPP1est2}
 \dbE^{\dbP^i} [\wh X^{t, \o^i}_{\t^i\wedge \ch^{t,\o^i}}]  - \dbE^{\dbP^{i,\o}} [\wh X^{t, \o}_{\tilde\t^i\wedge \ch^{t,\o}}]  
&\le&  C\rho_1(L\d).
\eea
Then by  Lemma \ref{lem-YUComega} and \reff{DPP1est1}, \reff{DPP1est2}  we have
\bea
\label{DPP1est3}
 \wh Y_t(\o)
 \;\le\;
 \wh Y_t (\o^i)+C\rho_1(L\d)
 \;\le\;
 \dbE^{\dbP^{i,\o}} [\wh X^{t, \o}_{\tilde\t^i\wedge \ch^{t,\o}}]  
 +\d +C\rho_1(L\d),\q \mbox{for all}~\o\in E_i.
 \eea

We next define:
\beaa
\tilde\t :=
  \1_{\{\t< t\}\cup\{\ch\le t\}}  \t 
 + \sum_{i\ge 1}\1_{E_i}\tilde \t^{i}(B^t), &\mbox{and then}& \{\t < t\}\cup\{ \ch\le t\} = \{\tilde\t < t\}\cup\{\ch\le t\}.
 \eeaa
Since $\t \le \ch$,  we see that $\{\t< t\}\cup\{\ch\le t\}=\{\t< t\}\cup\{\t=\ch = t\}$, and thus it is clear that  $\tilde\t \in \cT$.
  Moreover, we claim that there exists $\tilde \dbP\in \cP$ such that
 \bea
 \label{DPP1claim}
&\tilde \dbP = \dbP ~\mbox{on}~\cF_t ~~\mbox{and}~ \mbox{the r.c.p.d.}&\\
&(\tilde \dbP)^{t,\o} = \dbP^{i,\o}~\mbox{for}~\dbP\mbox{-a.e.}~\o\in E_i,  i\ge 1,\q (\tilde \dbP)^{t,\o} = \dbP^{t,\o}~\mbox{for}~\dbP\mbox{-a.e.}~\o\in \{\t< t\}\cup\{ \ch\le t\}.&\nonumber
 \eea
Then, by  \reff{DPP1est3} we have
 \bea
 \label{DPP1est4}
 \wh Y_t(\o) &\le& \dbE^{(\tilde\dbP)^{t,\o}}\big[\wh X^{t,\o}_{(\tilde \t \wedge \ch)^{t,\o}}\big] + \d + C\rho_1(L\d),\q   \dbP\mbox{-a.e.}~\o\in \{\t\ge t, \ch > t\},
 \eea
and therefore:
 \beaa
&& \dbE^{\dbP}
 \Big[\wh X_{\t}
      \1_{\{\t< t\}\cup\{\ch \le t\}} 
      +\wh Y_t \1_{\{\t\ge t, \ch > t\}}
 \Big]
 \\
 &\le& 
 \dbE^{\tilde\dbP}
 \Big[\wh X_{\tilde\t\wedge\ch}
      \1_{\{\t< t\}\cup\{\ch \le t\}} 
      +   \wh X_{\tilde\t\wedge \ch} \1_{\{\t\ge t, \ch > t\}}
 \Big]+\d +C\rho_1(L\d)\\
 &=& 
 \dbE^{\tilde\dbP}
 \Big[\wh X_{\tilde\t\wedge\ch} \Big]+\d +C\rho_1(L\d) 
 \;\le\; 
 \wh Y_0 + \d +C\rho_1(L\d),
 \eeaa
which implies \reff{DPP1est} by sending  $\d\to 0$.  Then the reverse inequality of \reff{ineq1-DPPdeterministic} follows from the  arbitrariness of $\dbP$ and $\t$.

 It remains to prove \reff{DPP1claim}. For any $\e>0$ and each $i\ge 1$, there exists a partition $\{E^i_j, j \ge 1\}$  of $E_i$ such that  $\|\o-\o'\|_t \le\e$ for any $\o, \o' \in E^i_j$. Fix an $\o^{ij}\in E^i_j$ for each $(i, j)$. By Property (P3) we may define $\tilde \dbP^\e \in \cP$ by:
 \beaa
 \tilde \dbP^\e := \dbP \otimes_t \Big[\sum_{i\ge 1} \sum_{j\ge 1}  \dbP^{i, \o^{ij}} \1_{E^i_j} +  \dbP\1_{ \{\t< t\}\cup\{ \ch\le t\}}\Big].
 \eeaa
By Property (P1), $\cP$ is weakly compact. Then $\tilde \dbP^{\e}$ has a weak limit $\tilde \dbP\in \cP$ as $\e\to 0$. 
We now show that $\tilde \dbP$ satisfies all the requirements in \reff{DPP1claim}.  Indeed, for any partition $0=s_0<\cds<s_{m} = t <s_{m+1} <\cds< s_{M} = t_\d < s_{M+1} <\cds<s_N=T$ and any bounded and uniformly continuous function $\f: \dbR^{N\times d} \to \dbR$, let $\xi := \f\big(B_{s_1} - B_{s_0},\cds, B_{s_N} - B_{s_{N-1}}\big)$. Then, denoting $\D s_k := s_{k+1}-s_k$, $\D\o_k:=\o_{s_k}-\o_{s_{k-1}}$, we see that 
\beaa
\dbE^{ \dbP^{i, \o}}[\xi^{t,\o}]
 =
 \eta^i_t(\o),\q \dbE^{ \dbP^{i, \o^{ij}}}[\xi^{t,\o}]
 =
 \eta^{i,j}_t(\o),
 \eeaa
  where:
 \beaa
 \eta^i_t(\o) 
 &:=&
 \dbE^{\dbP^i}\Big[\f\Big( (\D\o_k)_{1\le k \le m}, 
                          {\o^i_t - \o_t\over \d} 
                          (\D s_k)_{m+1\le k \le M},  
 (B_{s_{k}-\d} - B_{s_{k-1}-\d})_{M+1\le k\le N} \Big)\Big];\\
 \eta^{i,j}_t(\o) 
 &:=&
 \dbE^{\dbP^i}\Big[\f\Big( (\D\o_k)_{1\le k \le m}, 
                          {\o^i_t - \o^{ij}_t\over \d} 
                          (\D s_k)_{m+1\le k \le M},  
 (B_{s_{k}-\d} - B_{s_{k-1}-\d})_{M+1\le k\le N} \Big)\Big].
 \eeaa
Let $\rho$ denote the modulus of  continuity function of $\f$. Then
 \beaa
\Big| \dbE^{ \dbP^{i, \o^{ij}}}[\xi^{t,\o}] -\dbE^{ \dbP^{i, \o}}[\xi^{t,\o}]\Big|\le \rho(\e) &\mbox{for all}& \o\in E^i_j,
 \eeaa
% We emphasize that $\eta^i_t(\o)$ depends only on $\xi$, but does not depend on $n$ and $m_i$.  Then
 and thus
 \beaa
 &&\Big| \dbE^{\tilde \dbP^{\e}} [ \xi]   
 - \dbE^\dbP\big[\xi \1_{ \{\t< t\}\cup\{ \ch\le t\}}
                 +\sum_{i\ge 1} \eta^i_t  \1_{E_i}
            \big]
 \Big|\\
& =&\Big| \dbE^\dbP\Big[ \sum_{i,j\ge 1} 
                         \dbE^{\dbP^{i, \o^{ij}}}
                         [\xi^{t,\cd}] \1_{E^i_j}
                   \Big] 
          - \dbE^\dbP\Big[\sum_{i,j\ge 1}\eta^i_t\1_{E^i_j}
                     \Big]
     \Big|
 \\
 &\le& \dbE^\dbP\Big[ \sum_{i,j\ge 1} \big| \dbE^{ \dbP^{i, \o^{ij}}}[\xi^{t,\cd}] -\dbE^{ \dbP^{i, \cd}}[\xi^{t,\cd}]\big|\1_{E^i_j}\Big] 
\le  \dbE^\dbP\Big[ \sum_{i,j\ge 1}\rho(\e) \1_{E^i_j}\Big] 
 \;\le\; \rho(\e).
 \eeaa
By sending $\e\to 0$, we obtain
  $\dbE^{\tilde \dbP} [ \xi]   =\dbE^\dbP\big[\xi \1_{ \{\t< t\}\cup\{ \ch\le t\}}+\sum_{i\ge 1} \eta^i_t  \1_{E_i}\big]$, which proves \reff{DPP1claim} by the arbitrariness of $\xi$.
 \qed
 
\vspace{5mm}

We now prove the regularity in the $t$-variable. Recall the $\rho_2$ defined in \reff{rho12}.

\begin{lem}\label{lem-Ycontt}
%Let $X\in C^0_b(\L)$, $\textsc{h}\in\cH$ correspond to $O$ and $t_0$, and $\d>0$. Then, 
Let $0\le t_1 < \ch(\o^1)$, $0\le t_2 <\ch(\o^2)$, and $t_1 \le t_2$. Then  we have:
 \beaa
 |\wh Y_{t_1}(\o^1) - \wh Y_{t_2}(\o^2)| &\le& C\Big[1+{1\over d(\o^1_{t_1}, O^c)}\Big] \rho_2\Big( \dbf_\infty\big((t_1,\o^1), (t_2,\o^2)\big)\Big).
  \eeaa
\end{lem}

\proof Without loss of generality we assume  $t_1 < t_2$. Also, in view of the uniform continuity in $\o$ of Lemma \ref{lem-YUComega}, it suffices to prove the lemma in the case $\o^1=\o^2=\o$.

Denote $ \d :=  \dbf_\infty\big((t_1,\o), (t_2,\o)\big)$ and $\e := d(\o_{t_1}, O^c)$. For $\d \ge {1\over 8}$, we have
$ |\wh Y_{t_1}(\o) - \wh Y_{t_2}(\o)| \le 2C_0 \le C\e^{-1} \rho_2(\d)$. So we assume in the rest of this proof that $\d < {1\over 8}$. 

First, by Assumption \ref{assum-additional}, we may consider the measure $\dbP \in \cP_{t_1}$ such that $\a^{\dbP}_t := 0, \b^\dbP_t:=0$, $t\in [t_1, t_2]$. Then, by setting $\t := t_0$ in Lemma \ref{lem-DPP1}, we see that $\wh Y_{t_1}(\o) 
 \ge
 \cE_{t_1}[ \wh Y_{t_2}^{t_1,\o}]
 \ge
 \dbE^\dbP[ \wh Y_{t_2}^{t_1,\o}] 
 =
 \wh Y_{t_2}(\o_{\cd\wedge t_1}).
 $
 Note that $\ch(\o_{\cd\wedge t_1}) = t_0>t_2$.  Thus, by Lemma \ref{lem-YUComega}, 
 \bea
 \label{Ycontt1}
 \wh Y_{t_2}(\o) - \wh Y_{t_1}(\o) 
 &\le&
  C\rho_1\Big(\dbf_\infty\big((t_2, \o_{\cd\wedge t_1}), (t_2, \o)\big)\Big)\le C\rho_1(\d) \le C\rho_2(\d) .
 \eea

Next, for  arbitrary $\t \in \cT^{t_1}$, noting that $\wh X \le \wh Y$ we have
 \beaa
 I(\tau) &:=&\cE_{t_1}
   \Big[ \wh X^{t_1,\o}_{\t \wedge\textsc{h}^{t_1,\o}}
         \1_{\{\t \wedge \textsc{h}^{t_1,\o} < t_2\}}
         + \wh Y_{t_2}^{t_1,\o} 
           \1_{\{\t \wedge \textsc{h}^{t_1,\o} \ge  t_2\}}
   \Big]- \wh Y_{t_2}(\o)
 \\
  &=&\cE_{t_1}
   \Big[ \wh X^{t_1,\o}_{\t } \1_{\{\t <\textsc{h}^{t_1,\o}\wedge t_2 \}} + \wh X^{t_1,\o}_{\textsc{h}^{t_1,\o}}
         \1_{\{\textsc{h}^{t_1,\o} < t_2, \textsc{h}^{t_1,\o}\le \t \}}
         + \wh Y_{t_2}^{t_1,\o} 
           \1_{\{\t \wedge \textsc{h}^{t_1,\o} \ge  t_2\}}
   \Big]- \wh Y_{t_2}(\o)\\
 &\le& 
 \cE_{t_1}
 \Big[ \big(\wh X^{t_1,\o}_{\t} 
            - \wh X^{t_1,\o}_{\textsc{h}^{t_1,\o}\wedge t_2}
       \big)
       \1_{\{\t < \textsc{h}^{t_1,\o} \wedge t_2\}}
       + \wh Y^{t_1,\o}_{\textsc{h}^{t_1,\o}\wedge t_2}
 \Big]
 - \wh Y_{t_2}(\o)
 \\
  &\le& 
 \cE_{t_1}
 \Big[ \big(\wh X^{t_1,\o}_{\t} 
            - \wh X^{t_1,\o}_{\textsc{h}^{t_1,\o}\wedge t_2}
       \big)
       \1_{\{\t < \textsc{h}^{t_1,\o} \wedge t_2\}}
 \Big]  + \cE_{t_1}\Big[|\wh Y^{t_1,\o}_{t_2} 
 - \wh Y_{t_2}(\o)|
                            \1_{\{\textsc{h}^{t_1,\o} > t_2\}}
                       \Big]
 \\
 &&+C \;\cC_{t_1}\big[\textsc{h}^{t_1,\o} \le t_2\big].
 \eeaa
By \reff{USC2} and Lemma \ref{lem-YUComega} we have
 \beaa
 I(\tau)
 &\le&  
 \cE_{t_1}
 \Big[\rho_0\Big(\dbf_\infty ( (t_1, \o), (t_2, \o\otimes_{t_1} B^{t_1}))\Big)\Big] + C \cE_{t_1}
 \Big[\rho_1\big(\|\o- \o\otimes_{t_1} B^{t_1}\|_{t_2} \big)\Big]\\
&& 
 + C \cC_{t_1}\big[\|B^{t_1}\|_{t_2} 
                       \ge \e \big]\\
  &\le&  
 \cE_{t_1}
 \Big[\rho_0\big( \d + \|B^{t_1}\|_{t_2}\big)\Big]   
 +C  \cE_{t_1}
 \Big[\rho_1\big( \d + \|B^{t_1}\|_{t_2}\big)\Big]   + C\e^{-1} \cE_{t_1}\Big[ \|B^{t_1}\|_{t_2} \Big] \\
 &\le& C [1+\e^{-1}]  \cE_{t_1}
 \Big[\rho_1\big( \d + \|B^{t_1}\|_{t_2}\big)\Big] .
 \eeaa
Since $\d\le {1\over 8}$, following the proof of \reff{rho1est} we have
 \beaa
 I(\tau)
 &\le&  
 C[1+\e^{-1}]\Big[ \rho_1(\d^{1\over 3}) +  \d^{1\over 3}\Big]  
 \le  
 C[1+\e^{-1}] \rho_2(\d).
 \eeaa
By the arbitrariness of $\tau$ and the dynamic programming principle of Theorem \ref{thm-optimal-DP}, we obtain $\wh Y_{t_1}(\o) - \wh Y_{t_2}(\o) \le C\e^{-1}\rho_2(\d)$, and the proof is complete by \reff{Ycontt1}.
\qed

\vspace{5mm}
Applying  Lemmas \ref{lem-YUComega}, \ref{lem-DPP1}, and \ref{lem-Ycontt}, and following the same arguments as those of Theorem \ref{thm-T-optimal-DP},
we establish the dynamic programming principle in the present context.

\begin{thm}\label{thm-optimal-DP}
%Under Assumption \ref{assum-additional}, 
Let $t<\ch(\o)$ and $\t \in \cT^t$. Then 
 \beaa
 \wh Y_t(\o) 
& =&
 \sup_{\tilde\tau\in\cT^t}
 \cE_t
 \Big[\wh X^{t,\o}_{\tilde\t\wedge\textsc{h}^{t,\o}}
      \1_{\{\tilde\t\wedge\textsc{h}^{t,\o}<\t\}}
      +\wh Y^{t,\o}_{\t}
       \1_{\{\tilde\t\wedge\textsc{h}^{t,\o}\ge \t\}}\Big].
 \eeaa
Consequently, $\wh Y$ is a $\cE-$supermartingale on $[0,\textsc{h}]$.
\end{thm}

\ms
By Lemma \ref{lem-Ycontt}, $\wh Y$ is continuous for $t\in [0, \ch)$. Moreover, since $\hat Y$ is an $\cE$-supermartingale, we see that $\hat Y_{\ch-}$ exists. However, the following example shows that in general $\wh Y$ may be discontinuous at $\ch$.
\begin{eg}
\label{eg-discont}
{\rm Set $X_t(\o) := t$ and  let $\ch$ correspond to $O$ and $t_0$.  Clearly $\wh X = X$, $\wh Y_\ch = \ch$ and $\wh Y_t(\o) \le t_0$. However, for any $t < \ch(\o)$, set $\t := t_0$ and $\dbP\in \cP_t$ such that $\a^\dbP=0, \b^\dbP=0$, we see that $\wh Y_t(\o) \ge \dbE^\dbP\Big[X(\ch(\o\otimes_t B^t), \o\otimes_t B^t)\Big] = X(\ch(\o_{\cd\wedge t}), \o_{\cd\wedge t}) = \ch(\o_{\cd\wedge t}) = t_0$.
That is, $\wh Y_t(\o) = t_0$. Thus $\wh Y$ is discontinuous at $\ch$ whenever $\ch(\o) <t_0$.
\qed}
\end{eg}
This issue is crucial for our purpose, and we will discuss more in Subsection \ref{sect-prop-chcont} below.

\subsection{Continuous approximation of the hitting times}
\label{sect-hitting}

Similar to the proof of Theorem \ref{thm-T-optimal}, we need to apply some limiting arguments. % to the context of  Proposition \ref{prop-capacity} (iii). 
We therefore assume without loss of generality that $\wh Y_0 > \wh X_0$ and  introduce the stopping times: for any $m\ge 1$ and $n> (\wh Y_0 - \wh X_0)^{-1}$,
  \bea\label{Hn}
 \textsc{h}_m
 :=
 \inf\big\{t\ge 0: \dbf(\o_t, O^c) \le {1\over m}\big\}
 \wedge (t_0-\frac{1}{m}), \q \t_n := \inf\{t\ge 0: \wh Y_t - \wh X_t \le {1\over n}\}.
 \eea
 Here we abuse the notation slightly by using the same notation $\t_n$ as in \reff{T-tn}.
Our main task in this subsection is to build an approximation of $\textsc{h}_m$ and $\tau_n$ by continuous random variables. This will be obtained by a repeated use of Lemma \ref{lem-qscont}.

We start by a continuous approximation of the sequence $(\textsc{h}_m)_{m\ge 1}$ defined in \reff{Hn}.

\begin{lem}\label{lemhath}
For all $m\ge 2$:
\\
{\rm (i)}\quad $\textsc{h}_{m-1}(\o)
                \le 
                \textsc{h}_m(\o')
                \le
                \textsc{h}_{m+1}(\o)$, whenever $\|\o-\o'\|_{t_0}\le\frac{1}{m(m+1)}$,
\\
{\rm (ii)}\quad there exists an open subset $\O^m_0\subset\O$, and a uniformly continuous $\hat{\textsc{h}}_m$ such that
 \beaa\label{tildetaum}
 \cC\big[(\O^m_0)^c\big] <  2^{-m}  
 &\mbox{and}& 
 \textsc{h}_{m-1} - 2^{-m} 
 \le 
 \hat{\textsc{h}}_m
 \le 
 \textsc{h}_{m+1} + 2^{-m} 
 ~~\mbox{on}~~\O^m_0,
 \eeaa
{\rm (iii)}\quad there exist $\d_m>0$ such that $|\hat{\textsc{h}}_m(\o)-\hat{\textsc{h}}_m(\o')|\le 2^{-m}$ whenever $\|\o-\o'\|_{t_0}\le \d_m$, and:
 \beaa
 \cC\big[(\hat\O^m_0)^c\big] \le  2^{-m}
 &\mbox{where}&
 \hat\O^m_0
 :=
 \{\o\in\O^m_0:  \dbf(\o,[\O^m_0]^c) > \d_m\}.
 \eeaa
\end{lem}

\proof Notice that (ii) is a direct consequence of (i) obtained by applying Lemma \ref{lem-qscont} with $\e =  2^{-m}$. To prove (i), we observe that for $\|\o -\o'\|_{t_0} \le {1\over m(m+1)}$ and $t<\textsc{h}_m(\o')$, we have
 \beaa
 \dbf(\o_t, O^c) \ge \dbf(\o'_t, O^c) -  {1\over m(m+1)} 
      > {1\over m} - {1\over m(m+1)} = {1\over m+1}.
\eeaa
 This shows that $\textsc{h}_m(\o')\le\textsc{h}_{m+1}(\o)$ whenever $\|\o -\o'\|_{t_0} \le {1\over m(m+1)}$. Similarly, $\textsc{h}_{m-1}(\o) \le \textsc{h}_m(\o')$ whenever  $\|\o -\o'\|_{t_0} \le {1\over m(m-1)}$, and the inequality (i) follows. 

It remains to prove (iii). The first claim follows from the uniform continuity of $\hat{\textsc{h}}_m$. For each $\d>0$, define $h_\d: [0,\infty)\to [0,1]$ as follows: 
\bea
\label{hd}
\mbox{$h_\d(x) := 1$ for $x\le \d$, $h_\d(x) = 0$ for $x \ge 2\d$, and $h_\d$ is linear on $[\d, 2\d]$.}
\eea
 Then the map $\o\longmapsto\psi_\d(\o) := h_{\d}(d(\o, [\O^m_0]^c))$ is continuous, and $\psi_\d \downarrow \1_{[\O^m_0]^c}$ as $\d\downarrow 0$.  Applying Proposition \ref{prop-capacity} (ii) we have 
\beaa
\lim_{\d\to 0} \cE[\psi_\d] 
= 
\cE\big[ \1_{(\O^m_0)^c}\big] 
= 
\cC\big[(\O^m_0)^c\big]< 2^{-m}.
\eeaa
By definition of $\hat\O^m_0$, notice that $\1_{(\hat\O^m_0)^c}\le\psi_{\delta_m}$. Then $\cC\big[(\hat\O^m_0)^c\big] \le  \cE[\psi_{\d_m}]$, and (iii) holds true for sufficiently small $\d_m$.
\qed

\vspace{5mm}

We next derive a continuous approximation of the sequences
 \bea\label{taunm}
 \tau^m_n
 &:=& \tau_n\wedge\hat{\textsc{h}}_m,
% &n,m\ge 1.& 
 \eea
where $\tau_n$ and $\hat{\textsc{h}}_m$ are defined in \reff{Hn} and Lemma \ref{lemhath} (ii), respectively.

\begin{lem}\label{lemtaumn}
For all $m\ge 2$, $n > (\wh Y_0-\wh X_0)^{-1}$, there exists an open subset $\O_n^m\subset\O$ and a uniformly continuous map $\hat\tau^m_n$ such that
 $$
 \tau^m_{n-1}-2^{1-m}-2^{-n}
 \le
 \hat\tau^m_n
 \le
 \tau^m_{n+1}+2^{1-m}+2^{-n}
 ~\mbox{on}~
 \hat\O_0^m\cap\O_n^m,
 ~\mbox{and}~
 \cC\big[(\O^m_n)^c\big] 
 \le 
 2^{-n}.
 $$
\end{lem}

\proof
Fix $m$, and recall the modulus of continuity $\rho_1$ introduced in \reff{rho12}. For each $n$, let $0< \d^m_n < \d^m$ such that  $(\rho_0+C\rho_1)(\d^m_n) \le {1\over n(n+1)}$, where $C$ is the constant in Lemma  \ref{lem-YUComega} .  We shall prove
\bea
\label{hatmnest}
&(\t_{n-1}\wedge\hat{\textsc{h}}_m)(\o)-2^{1-m} 
\le 
(\t_n\wedge\hat{\textsc{h}}_m)(\o') 
\le 
(\t_{n+1}\wedge\hat{\textsc{h}}_m)(\o)+2^{1-m}&
\\
&\mbox{whenever}~ \o\in \hat\O^m_0, 
 ~~\|\o -\o'\|_{t_0} \le \d^m_n.&
 \nonumber
\eea
Then the required statement follows from Lemma \ref{lem-qscont} with $\e = 2^{-n}$.

We shall prove only the right inequality of \reff{hatmnest}. The left one can be proved similarly.  Let $\o, \o'$ be as in \reff{hatmnest}. First,  by Lemma \ref{lemhath} (iii) we have
\bea
\label{hatmnest1}
\o' \in \O^m_0 &\mbox{and}& \hat{\textsc{h}}_m(\o') 
\le
\hat{\textsc{h}}_m(\o)+2^{-m}
\eea
We now prove the right inequality of \reff{hatmnest} in three cases.
\\
{\it Case 1.} if $\t_{n+1}(\o)  \ge \hat\ch_m(\o') - 2^{-m}$, then $\hat{\textsc{h}}_m(\o') 
\le
(\t_{n+1}\wedge \hat{\textsc{h}}_m)(\o)+2^{-m}$ and thus the result is true.
\\
{\it Case 2.} If $\t_{n+1}(\o) =\ch(\o)$, then by Lemma \ref{lemhath} (ii) we have $\hat\ch_m(\o) \le \ch_{m+1}(\o) + 2^{-m} \le \t_{n+1}(\o) + 2^{-m}$, and thus $\hat\ch_m(\o') \le \hat{\textsc{h}}_m(\o)+2^{-m} \le \t_{n+1}(\o) + 2^{1-m}$. This, together with \reff{hatmnest1}, proves the desired inequality.
\\
{\it Case 3.} We now assume $\t_{n+1}(\o)  < \hat\ch_m(\o') - 2^{-m}$ and  $\t_{n+1}(\o) <\ch(\o)$. By Lemma \ref{lemhath} (ii) we have $\t_{n+1}(\o) <\textsc{h}_{m+1}(\o')$, and thus $\t_{n+1}(\o)< \ch(\o') $. Then it follows from Lemma \ref{lem-YUComega} that
\beaa
(Y-X)_{\t_{n+1}(\o)}(\o') \le (Y-X)_{\t_{n+1}(\o)}(\o) + (\rho_0+C\rho_1)(\d^m_n) \le {1\over n+1} + {1\over n(n+1)} = {1\over n}.
\eeaa
That is, $\t_n(\o') \le \t_{n+1}(\o)$. This, together with \reff{hatmnest1}, proves the desired inequality.
\qed

\vspace{5mm}

For our final approximation result, we introduce the notations:
 \bea\label{taun*}
 \bar\tau_n:=\tau_n\wedge\textsc{h}_n,
 ~~\underline\theta^*_n:=\hat\tau_{n-1}^{n-1}-2^{3-n},
 ~~\overline\theta^*_n:=\hat\tau_{n+1}^{n+1}+2^{1-n},
 \eea
and
 \bea\label{Omegan*}
 \Omega^*_n:=\hat\O^{n-1}_0\cap\O_{n-1}^{n-1}\cap
             \hat\O^{n+1}_0\cap\O_{n+1}^{n+1}.
 \eea

\begin{lem}
\label{lem-Omegan*}
For all $n\ge (\wh Y_0-\wh X_0)^{-1}\vee 2$, $\underline\theta^*_n,\overline\theta^*_n$ are uniformly continuous, and $\underline\theta^*_n\le\bar\tau_n\le\overline\theta^*_n$ on $\O^*_n$.
\end{lem}

\proof This is a direct combination of Lemmas \ref{lemhath} and \ref{lemtaumn}.

\subsection{Proof of Theorem \ref{thm-optimal} }
%\label{sect-Proof}
%\setcounter{equation}{0}
We first prove the $\cE$-martingale property under an additional condition.

\begin{lem}
\label{lem-mg}
Let $\t \in \cT$ such that $\t \le \t^*$ and $\cE[Y_{\t-} ] = \cE[Y_\t]$ (in particular if $\t < \ch$). Then $\wh Y$ is an $\cE$-martingale on $[0, \t]$.
\end{lem}
\proof  If $\wh Y_0=\wh X_0$, then $\wh\t^*=0$ and obviously the statement is true. We then assume $\wh Y_0>\wh X_0$, and prove the lemma in several steps.
\\
{\it Step 1}\quad Let $n$ be sufficiently large so that ${1\over n} < \wh Y_0-\wh X_0$.  Follow the same arguments as that of Lemma \ref{lem-T-DPn} , one can easily prove:
 \bea\label{ELmartn}
 \wh Y
 &\mbox{is an}& 
 \cE-\mbox{martingale on}
 ~~[0,\t_n].
 \eea 
\no {\it Step 2}\quad Recall the sequence of stopping times $(\bar\tau_n)_{n\ge 1}$ introduced in \reff{taun*}. By Step 1 we have $\wh Y_0 = \cE[\wh Y_{\bar\t_n}]$. Then for any $\e>0$, there exists $\dbP_n\in \cP$ such that $\wh Y_0 - \e < \dbE^{\dbP_n}[\wh Y_{\bar\t_n}]$. Since $\cP$ is weakly compact, there exists subsequence $\{n_j\}$ and $\dbP^*\in \cP$ such that $\dbP_{n_j}$ converges weakly to $\dbP^*$. Now for any $n$ and $n_j\ge n$, since $Y$ is a supermartingale under each $\dbP_{n_j}$ and $(\bar\t_n)_{n\ge 1}$ is increasing, we have 
 \bea
 \label{Y0e}
 \wh Y_0 - \e 
 &<& 
 \dbE^{\dbP_{n_j}}\big[\wh Y_{\bar\t_{n_j}}\big] 
 \;\le\; 
 \dbE^{\dbP_{n_j}}\big[\wh Y_{\bar\t_{n}}\big].
 \eea
Our next objective is to send $j\nearrow\infty$, for fixed $n$, and use the weak convergence of $\dbP^{n_j}$ towards $\dbP^*$. To do this, we need to approximate $\wh Y_{\bar\t_{n}}$ with continuous random variables. Denote 
 \bea\label{psimn}
 \psi_n(\o) 
 := 
 h_n\Big(\inf_{0\le t\le \overline\th^*_n(\o)}
         \dbf(\o_t, O^c) \Big)
 ~\mbox{with}~
 h_n(x):=1\wedge[(n+3)(n+4)x-(n+3)]^+.
 \eea
 Then $\psi_n$ is continuous in $\o$, and 
\bea
\label{psinthn}
 &\{\psi_n >0\} 
  \subset 
  \big\{\inf_{0\le t\le \overline\th^*_n(\o)}
        \dbf(\o_t, O^c) > {1\over n+4}
  \big\} 
  \subset  
  \{\overline\th^*_n < \ch_{n+4}\}.
\eea
 In particular, this implies that $\wh Y_{\underline\th^*_n} \psi_n$ and  $\wh Y_{\overline\th^*_n} \psi_n$ are continuous in $\o$. We now decompose the right hand-side term of \reff{Y0e} into:
 \beaa
% \label{Y0e2}
 \wh Y_0 - \e 
 &\le& 
 \dbE^{\dbP_{n_j}}
 \Big[\big[\wh Y_{\underline\th_{n}^*}  
           +(\wh Y_{\bar\t_{n}}-\wh Y_{\underline\th_{n}^*})
            \1_{\O^*_n} 
      \big]\big(\psi_n +(1-\psi_n)\big)  
      +(\wh Y_{\bar\t_{n}}-\wh Y_{\underline\th_{n}^*})
       \1_{(\O^*_n)^c}
 \big].
  \eeaa
Note that $\underline\th_n^* \le \bar\t_n \le \overline\th_n^*$  on $\O^*_n$. Then
 \beaa
 \wh Y_0 - \e 
 &\le& 
 \dbE^{\dbP_{n_j}}
 \Big[\Big(\wh Y_{\underline\th_{n}^*}  
                 + \sup_{\underline\th_n^* 
                         \le t \le \overline\th_n^*}
                   (\wh Y_t-\wh Y_{\underline\th_{n}^*})
             \Big)\psi_n
 \Big] 
 + C\cC [\psi_n<1] 
 + C\cC\big[(\O^*_n)^c\big].
 \eeaa
 Send $j\to\infty$, we obtain 
  \bea
 \label{Y0e2}
 \wh Y_0 - \e 
 \le
 \dbE^{\dbP^*}
 \Big[\psi_n\wh Y_{\underline\th_{n}^*}\Big]
 + \dbE^{\dbP^*}
 \Big[\psi_n\sup_{\underline\th_n^* \le t \le \overline\th_n^*}
            (\wh Y_t-\wh Y_{\underline\th_{n}^*})
 \Big] 
 + C\cC [\psi_n<1] 
 + C \cC[(\O^*_n)^c].
 %\nonumber\\ 
 %&=:& A_n + B_n + C_n + D_n.
  \eea
%where
 %\beaa
 %A_{n,j}:=\dbE^{\dbP_{n_j}}\big[Y_{\underline\th_{n}^*}\big],
% &B_{n,j}:=\dbE^{\dbP_{n_j}}
  %         \big[(Y_{\t_{n}^*}-Y_{\underline\th_{n}^*})
 %               \1_{\O^*_n}\big],&
% C_n:=2\|Y\|_\infty \;\sup_{\dbP\in{\cP}^L}
 %                     \dbP\big[(\O^*_n)^c\big].
 %\eeaa
\no {\it Step 3.} In this step we show that
  \bea\label{BCD}
 \lim_{n\to\infty} \dbE^{\dbP^*}
 \Big[\psi_n\sup_{\underline\th_n^* \le t \le \overline\th_n^*}
            (\wh Y_t-\wh Y_{\underline\th_{n}^*})
 \Big] 
 = \lim_{n\to\infty} \cC [\psi_n<1] 
 =\lim_{n\to\infty} \cC[(\O^*_n)^c] =0.
 \eea 
 (i) First, by the definition of $\O_n^*$ in \reff{Omegan*} together with Lemmas \ref{lemhath} (iii) and \ref{lemtaumn}, it follows that $\cC\big[(\O^*_n)^c\big] \le C 2^{-n}\longrightarrow 0$ as $n\to\infty$.
\\
(ii) Next, notice that 
 \beaa
 &
 \{\psi_n < 1\} 
 \;=\;  
 \big\{\inf_{0\le t\le \overline\th^*_n(\o)}
       \dbf(\o_t, O^c) 
       < {1\over n+3}
 \big\} 
 \;\subset\;
 \{\overline\th_n^* > \ch_{n+3}\} .&
 \eeaa
Moreover,  by \reff{taun*} and Lemma \ref{lem-Omegan*},
 \beaa
 \overline\th_n^* 
 =
 \hat \t^{n+1}_{n+1} +2^{1-n} 
 = 
 \underline\th_{n+2}^* + 2^{2-n} 
 \le 
 \bar\t_{n+2} + 2^{2-n} 
 \le 
 \ch_{n+2} + 2^{2-n}, 
 &\mbox{on}& \O^*_{n+2}.
\eeaa
Then
 $$
 \begin{array}{rl}
 \{\psi_n < 1\}&
 \subset\;
 (\O^*_{n+2})^c \cup \{\ch_{n+3} < \ch_{n+2} + 2^{2-n}\}
 \\
 &
 \subset\; 
 (\O^*_{n+2})^c 
 \cup 
 \Big\{\sup_{\ch_{n+2} \le t\le \ch_{n+2}+2^{2-n}}    
       |B_t-B_{\ch_{n+2}}| 
       \ge 
       {1\over (n+2)(n+3)}
 \Big\} .
 \end{array}
 $$
Then one can easily see that $ \cC[\psi_n<1]\to 0$, as $n\to \infty$.
\\
(iii) Finally, it is clear that $\underline \th_n^*\to \wh\t^*$, $\overline\th^*_n \to \wh\t^*$. Recall that $\wh Y_{\wh \t^*-}$ exists. By \reff{psinthn}, we see that $\psi_n \sup_{\underline\th_n^* \le t \le \overline\th_n^*}  (\wh Y_t-\wh Y_{\underline\th_{n}^*})\to 0$, $\dbP^*$-a.s. as $n\to\infty$. Then by applying the dominated convergence theorem under $\dbP^*$ we obtain the first convergence in \reff{BCD}.
\\
\no {\it Step 4.} By the dominated convergence theorem under $\dbP^*$ we obtain
$
 \lim_{n\to\infty} \dbE^{\dbP^*}
 [\psi_n\wh Y_{\underline\th_{n}^*}] = \dbE^{\dbP^*}[\wh Y_{\wh\t^*-}] .
$ 
This, together with \reff{Y0e2} and \reff{BCD}, implies that
\beaa
\wh Y_0 \le \dbE^{\dbP^*}[\wh Y_{\wh\t^*-}] + \e.
\eeaa
Note that $\wh Y$ is an $\dbP^*$-supermartingale and $\t\le \wh\t^*$, then
\beaa
\wh Y_0 \le \dbE^{\dbP^*}[\wh Y_{\t-}] + \e.
\eeaa
Since $\e$ is arbitrary, we obtain  $\wh Y_0 \le \cE[\wh Y_{\t-}] $, and thus by the assumption  $ \cE[\wh Y_{\t-}] =  \cE[\wh Y_{\t}]$ we have $\wh Y_0 \le \cE[\wh Y_{\t}] $. This, together with the fact that $\wh Y$ is a $\cE$-supermartingale, implies that 
\bea
\label{mgtau*}
\wh Y_0 = \cE[\wh Y_{\t}] .
\eea
Similarly, one can prove $\wh Y_t(\o) = \cE_t[\wh Y^{t,\o}_{\t^{t,\o}}]$ for $t<\t(\o)$, and thus $\wh Y_{.\wedge\tau}$ is a $\cE$-martingale.
\qed

\vspace{5mm}

In light of Lemma \ref{lem-mg}, the following result is obviously  important for us.

\begin{prop}
\label{prop-chcont}
It holds that $\cE[\wh Y_{\wh\t^*-}] = \cE[\wh Y_{\wh\t^*}]$.
\end{prop}

We recall again that $\wh Y_{\wh\t^*-}=\wh Y_{\wh\t^*}$ whenever $\wh\t^* <\ch$. So the only possible discontinuity is at $\ch$. The proof of Proposition \ref{prop-chcont} is reported in Subsection \ref{sect-prop-chcont} below. Let us first show how it allows to complete the

\vspace{5mm}

\no {\bf Proof of Theorem \ref{thm-optimal}}
By Lemma \ref{lem-mg} and  Proposition \ref{prop-chcont},  $\wh Y$ is an $\cE$-martingale on $[0, \wh\t^*]$. Moreover, since $\wh X_{\wh\t^*} = \wh Y_{\wh\t^*}$, then $\wh Y_0 = \cE[\wh X_{\wh\t^*}]$ and thus $\wh\t^*$ is an optimal stopping time.
\qed

\subsection{$\cE-$Continuity of $\wh Y$ at the random maturity}
\label{sect-prop-chcont}

This subsection is dedicated to the proof of Proposition \ref{prop-chcont}.  We first reformulate some pathwise properties established in previous subsections. For that purpose, we introduce the following additional notation:  for any  $\dbP\in\cP$, $\t\in\cT$, and $ E\in \cF_\t$
 \bea
 \label{cPtau}
 \cP(\dbP,\t, E)
 := \Big\{ \dbP' \in \cP: \dbP' =  \dbP \otimes_\t \big[ \dbP' \1_E + \dbP \1_{E^c}\big]\Big\},\q \cP(\dbP,\t) := \cP(\dbP,\t, \O).
 \eea
 That is, $\dbP'\in \cP(\dbP,\t, E)$ means $\dbP' = \dbP$ on $\cF_\t$ and $(\dbP')^{\t, \o} = \dbP^{\t, \o}$ for $\dbP$-a.e. $\o\in E^c$. 

The first result corresponds to Theorem \ref{thm-optimal-DP}.

\begin{lem}
\label{lem-ch1}
Let $\dbP\in \cP$, $\t_1, \t_2 \in \cT$,  and $E\in \cF_{\t_1}$. Assume $\t_1 \le \t_2 \le \ch$, and $\t_1 <\ch$ on $E$. Then for any $\e>0$, there exist $\dbP_\e \in \cP(\dbP, \t_1, E)$ and $\t_\e\in \cT$ with values in $[\t_1, \t_2]$,  s.t. 
\beaa
\dbE^{\dbP}\Big[ \wh Y_{\t_1}\1_E\Big] \le \dbE^{\dbP_\e}\Big[\big[\wh X_{\t_\e}\1_{\{\t_\e < \t_2\}} + \wh Y_{\t_2} \1_{\{\t_\e = \t_2\}}\big] \1_E\Big] + \e.
\eeaa
\end{lem}

\proof  Let $\t^n_1$ be a sequence of stopping times such that $\t^n_1 \downarrow \t$ and each $\t^n_1$ takes only finitely many values. Applying Lemma \ref{lem-Ycontt} together with the dominated convergence Theorem under $\dbP$, we see that
$
%\lim_{n\to\infty}  \dbP(\t^n_1 \ge \ch) = 0 &\mbox{and thus}&
 \lim_{n\to\infty} \dbE^\dbP\Big[|\wh Y_{\t^n_1\wedge \t_2}- \wh Y_{\t_1}|\Big]  =0.
$
Fix $n$  such that
\bea
\label{tau1n}
\dbE^\dbP\Big[|\wh Y_{\t^n_1\wedge \t_2}- \wh Y_{\t_1}|\Big]  \le  {\e\over 2}. % &\mbox{and}& C_0 \dbP(\t_2 < \t^n_1) \le {\d\over 2}.
\eea
Assume $\t^n_1$ takes values $\{t_i, i=1,\cds, m\}$, and for each $i$,  denote $E_i := E\cap \{\t^n_1 = t_i <\t_2\}\in \cF_{t_i}$. By  \reff{DPP1est4}, there exists $\tilde\t_i\in \cT$ and $\tilde\dbP_i \in \cP(\dbP, t_i)$ such that $\tilde\t_i \ge t_i$ on $E_i$  and 
\bea
\label{ch1-est}
 \wh Y_{t_i} \le \dbE^{\tilde\dbP_i}_{t_i}\Big[\wh X_{\tilde \t_i\wedge \ch}\Big] + {\e\over 2},\q \dbP\mbox{-a.s. on}~E_i.
 \eea
%$E^j_i:= \{\ch_j \wedge \t_2\le \t_1^n = t_i < \ch_{j+1}\wedge \t_2\} \in \cF_{t_i}$, $j\ge 0$, where we take the convention that $\ch_0 :=0$. By Lemma  \ref{lem-YUComega} $\wh Y_{t_i}$ is uniformly continuous in $\o$ on $E^j_i$. Applying the dynamic programming principle of Theorem \ref{thm-optimal-DP} with $\t = \t_2$, and following standard arguments, there exist $t_i \le \t^j_i \le \t_2$ and $\dbP^j_i \in \cP(\dbP, t_i)$ such that 
 %\beaa
 %\wh Y_{t_i} \le \dbE^{\dbP^m_i}_{t_i}\Big[\wh X_{\t^m_i}\1_{\{\t^m_i < \t_2\}} + \wh Y_{\t_2}\1_{\{\t^m_i = \t_2\}}\Big] + {\e\over 2},\q \dbP\mbox{-a.s. on}~E^m_i.
 %\eeaa
 Here $\dbE^{\tilde\dbP_i}_{t_i}[\cd] := \dbE^{\tilde\dbP_i}[\cd | \cF_{t_i}]$ denotes the conditional expectation. Define 
 \bea
 \label{ch1-tildeP}
 \tilde\t 
 := 
 \t_2\1_{E^c\cup \{\t_2 \le \t_1^n\}} 
 + \sum_{i=1}^m  \tilde \t_i \1_{E_i}, &&
 \tilde\dbP
 := \dbP \1_{E^c\cup\{ \t_2 \le \t_1^n\}} 
 + \sum_{i=1}^m  \tilde \dbP_i \1_{E_i}.
% \dbP\big[E\cap \{\t_2 \le \t_1^n\}\big] 
 %+ \sum_i \sum_{m\ge 0} \dbP^m_i\big[E\cap E^m_i\big],
 %\q\mbox{for any}~E\in \cF_T.
 \eea
Then one can check straightforwardly that 
\bea
\label{ch1-tildet}
\tilde\t\in \cT &\mbox{and}& \tilde\t \ge \t_2 \wedge \t^n_1;
\eea
and  $\tilde\dbP \in \cP(\dbP,\t_2\wedge\t_1^n, E)\subset\cP(\dbP,\t_1, E)$. Moreover,  by \reff{ch1-est} and \reff{ch1-tildeP},
\beaa
 \dbE^{\tilde\dbP}\big[\wh Y_{\t_2\wedge \t_1^n}\1_{E}\big]  &=&  \dbE^{\tilde\dbP}\Big[\big[\wh Y_{\t_2} \1_{\{\t_2 \le \t^n_1\}} +\sum_{i=1}^m \wh Y_{t_i} \1_{E_i}\big]\1_{E}\Big]\\ 
 &\le&    \dbE^{\tilde\dbP}\Big[ \big[\wh Y_{\t_2} \1_{\{\t_2 \le \t^n_1\}} +(\wh X_{\tilde \t \wedge \ch} +{\e\over 2})\1_{\{\t^n_1<\t_2\}}\big]\1_{E}\Big].
\eeaa
This,  together with \reff{tau1n} and \reff{ch1-tildet} , leads to
\beaa
&& \dbE^{\tilde\dbP}
 \Big[\big(\wh Y_{\t_1} 
           -\wh X_{\tilde\t}\1_{\{\tilde\t < \t_2\}} 
           -\wh Y_{\t_2} \1_{\{\tilde\t\ge \t_2\}} 
      \big)\1_{E}
 \Big] \\
 &\le&  \e+ \dbE^{\tilde\dbP}
 \Big[\big(\wh Y_{\t_2} \1_{\{\t_2 \le \t^n_1\}} + \wh X_{\tilde \t \wedge \ch} \1_{\{\t^n_1<\t_2\}}
           -\wh X_{\tilde\t}\1_{\{\tilde\t < \t_2\}} 
           -\wh Y_{\t_2} \1_{\{\tilde\t \ge \t_2\}} 
      \big)\1_{E}
 \Big] \\
 &=&\e +  \dbE^{\tilde\dbP}
 \Big[\big( \wh X_{\tilde \t \wedge \ch} - \wh Y_{\t_2}\big) \1_{\{\t^n_1<\t_2\le \tilde\t\}}\1_{E}
 \Big] \\
 &=& \e +  \dbE^{\tilde\dbP}
 \Big[\big( \dbE^{\tilde\dbP}_{\t_2}[\wh X_{\tilde \t \wedge \ch}] - \wh Y_{\t_2}\big) \1_{\{\t^n_1<\t_2\le \tilde\t\}}\1_{E}
 \Big] 
 \le  \e,
 \eeaa
 where the last inequality follows from the definition of $\wh Y$.  Then, by setting $\t_\e := \tilde\t \wedge \t_2$ we prove the result.
\qed

Next result corresponds to Lemma \ref{lem-mg}.

\begin{lem}
\label{lem-ch2}
Let $\dbP\in \cP$, $\t\in \cT$, and $E\in \cF_\t$  such that $\t \le \wh\t^*$ on $E$. Then for all $\e>0$:
\beaa
\dbE^\dbP\big[ \1_E\wh Y_\t \big]
\le \dbE^{\dbP_\e}\big[\1_E\wh Y_{\wh \t^*-}\big] + \e
&\mbox{for some}&
\dbP_\e \in \cP(\dbP, \t, E).
\eeaa
\end{lem}
\proof  %Note that $\wh Y_{\wh\t^*} \le \wh Y_{\wh\t^*-}$ and $\dbP_\e = \dbP$ on $\{\t=\wh\t^*\}$, then
%\beaa
%\dbE^\dbP[ \wh Y_\t ] - \dbE^{\dbP_\e}[\wh Y_{\wh \t^*-}] = \dbE^\dbP[ \wh Y_\t \1_{\{\t<\wh\t^*\}}] - \dbE^{\dbP_\e}[\wh Y_{\wh \t^*-}\1_{\{\t<\wh\t^*\}}].
%\eeaa 
We proceed in three steps.
\\
{\it Step 1.} We first assume $\t = t<\wh\t^*$ on $E$.  We shall prove the result following the arguments in Lemma \ref{lem-mg}. Recall the notations in Subsection \ref{sect-hitting} and the $\psi_n$ defined in \reff{psimn}, and let $\rho_n$ denote the modulus of continuity functions of $\underline \th^*_n$, $\overline \th^*_n$, and $\psi_n$. 
%Note that $\bar\t_n \uparrow \wh\t^*$, so there exists $n_1$ such that
%\bea
%\label{ch2-est1}
%2C_0 \dbP(\bar\t_n \le \t)\le {\e\over 3} &\mbox{for all}& n\ge n_1.
%\eea

 Denote $\bar \t_n :=0$ for $n \le (\wh Y_0-\wh X_0)^{-1}$. For any $n$ and $\d>0$, let $\{E^{n,\d}_i, i\ge 1\}\subset \cF_t$ be a partition of $E\cap \{\bar\t_{n-1}\le t < \bar\t_n\}$ such that $\|\o-\o'\|_{t} \le \d$ for any $\o, \o'\in E^{n,\d}_i$.  For each $(n,i)$, fix $\o^{n,i} := \o^{n, \d, i}\in E^{n,\d}_i$. By Lemma \ref{lem-mg}, $\wh Y\1_{E^{n,\d}_i}$ is an $\cE$-martingale on $[t, \bar \t_n]$. Then $\wh Y_{t}(\o^{n,i}) = \cE_{t}[\wh Y^{t,\o^{n,i}}_{\bar\t_n^{t,\o^{n,i}}}]$, and thus there exists $\dbP^{n,\d}_i \in \cP_{t}$ such that 
\bea
\label{Pndi}
\wh Y_{t}(\o^{n,i}) \le \dbE^{\dbP^{n,\d}_i}\Big[\wh Y^{t,\o^{n,i}}_{\bar\t_n^{t,\o^{n,i}}}\Big] + \e.
\eea
Note that $\cup_{m=1}^n \cup_{i\ge 1} E^{m,\d}_i =E\cap \{t < \bar\t_n\}$. Set
\bea
\label{Pnd}
 \dbP^{n,\d} &:=& \dbP \otimes_t \Big[\sum_{m=1}^n \sum_{i\ge 1}\dbP^{m,\d}_i \1_{E^{m,\d}_i} + \dbP \1_{E^c\cup \{t \ge \bar\t_n\}}\Big] \in \cP(\dbP, t, E).
 \eea
 Recall the $h_\d$ defined by \reff{hd}.
 We claim that, for any $N\ge n$,
 \bea
 \label{ch2-claim}
&&\dbE^\dbP[ \wh Y_t \1_E] - \dbE^{\dbP^{N,\d}}[\wh Y_{t \vee \underline\th^*_n}\psi_n \1_E]\nonumber\\
&\le&  Cn\cE\Big[\rho_2\Big(\d + \rho_n(\d) + 2 \eta_n(\d)\Big) \Big]+ C\rho_n(\d)+\e + C 2^{-n} + C\cC(\psi_n < 1)\nonumber\\
&&+ 2 \dbE^{\dbP^{N,\d}}\Big[\sup_{\underline\th^*_n \le s \le \overline\th^*_n}|\wh Y_s- \wh Y_{\underline\th^*_n} |  \psi_n\1_E\Big] +C  \cE\Big[h_\d\big(d\big(\o,  (\O^*_n)^c\big)\big) \Big],\\
&\mbox{where}& \eta_n(\d):= \sup_{t\le s_1< s_2\le t_0, s_2 - s_1 \le \rho_n(\d)}|B^t_{s_1} - B^t_{s_2} |.\nonumber
 \eea
 Moreover, one can easily find $\cF_t$-measurable continuous random variables $\f_k$ such that  $|\f_k|\le 1$ and $\lim_{k\to\infty}\dbE^\dbP[|\1_E - \f_k|]=0$. Then
\beaa
&&\dbE^\dbP[ \wh Y_t \1_E] - \dbE^{\dbP^{N,\d}}[\wh Y_{t \vee \underline\th^*_n}\psi_n \f_k] \\
&\le& Cn\cE\Big[\rho_2\Big(\d + \rho_n(\d) + 2 \eta_n(\d)\Big) \Big]+ C\rho_n(\d)+\e + C 2^{-n} + C\cC(\psi_n < 1) \\
&&+ C \dbE^{\dbP^{N,\d}}\Big[\sup_{\underline\th^*_n \le s \le \overline\th^*_n}|\wh Y_s- \wh Y_{\underline\th^*_n} |  \psi_n\f_k\Big] +C\cE\Big[h_\d\big(d\big(\o,  (\O^*_n)^c\big)\big) \Big]+ C\dbE^\dbP[|\1_E - \f_k|].
\eeaa
Send $\d\to 0$. First note that $[\d + \rho_n(\d) + 2 \eta_n(\d)] \downarrow 0$ and $h_\d\downarrow \1_{\{0\}}$, then by Proposition \ref{prop-capacity} (ii) we have
\beaa
&&\lim_{\d\to 0} \cE\Big[\rho_2\Big(\d + \rho_n(\d) + 2 \eta_n(\d)\Big) \Big] = 0;\\
&&\lim_{\d\to 0} \cE\Big[h_\d\big(d\big(\o,  (\O^*_n)^c\big)\big) \Big] = \cC\Big[d\big(\o,  (\O^*_n)^c\big) = 0\Big] = \cC[(\O^*_n)^c] \le C 2^{-n}.
\eeaa
Moreover,
for each $N$,  by the weak compactness assumption (P1) we see that  $\dbP^{N,\d}$ has a weak limit $\dbP^N\in \cP$.  It is straightforward to check that $\dbP^N \in   \cP(\dbP, t, E)$. Note that the random variables $\wh Y_{t \vee \underline\th^*_n}\psi_n \f_k$ and $\sup_{\underline\th^*_n \le s \le \overline\th^*_n}|\wh Y_s- \wh Y_{\underline\th^*_n} |  \psi_n\f_k$ are continuous. Then
\beaa
&&\dbE^\dbP[ \wh Y_t \1_E] - \dbE^{\dbP^{N}}[\wh Y_{t \vee \underline\th^*_n}\psi_n \f_k] \\
&\le& \e + C 2^{-n} + C\cC(\psi_n < 1) + C \dbE^{\dbP^{N}}\Big[\sup_{\underline\th^*_n \le s \le \overline\th^*_n}|\wh Y_s- \wh Y_{\underline\th^*_n} |  \psi_n\f_k\Big] + C\dbE^\dbP[|\1_E - \f_k|].
\eeaa
Again by the weak compactness assumption (P1),    $\dbP^{N}$ has a weak limit $\dbP^*\in \cP(\dbP, t, E)$ as  $N\to\infty$. Now send $N\to \infty$, by the continuity of the random variables we obtain
\beaa
&&\dbE^\dbP[ \wh Y_t \1_E] - \dbE^{\dbP^*}[\wh Y_{t \vee \underline\th^*_n}\psi_n \f_k] \\
&\le& \e + C 2^{-n} + C\cC(\psi_n < 1) + C\dbE^{\dbP^{*}}\Big[\sup_{\underline\th^*_n \le s \le \overline\th^*_n}|\wh Y_s- \wh Y_{\underline\th^*_n} |  \psi_n\f_k\Big] + C\dbE^\dbP[|\1_E - \f_k|].
\eeaa
Send $k\to\infty$ and recall that $\dbP^* = \dbP$ on $\cF_t$, we have
\beaa
&&\dbE^\dbP[ \wh Y_t \1_E] - \dbE^{\dbP^*}[\wh Y_{t \vee \underline\th^*_n}\psi_n 1_E] \\
&\le& \e + C 2^{-n} + C\cC(\psi_n < 1) + 2 \dbE^{\dbP^{*}}\Big[\sup_{\underline\th^*_n \le s \le \overline\th^*_n}|\wh Y_s- \wh Y_{\underline\th^*_n} |  \psi_n\1_E\Big] .
\eeaa
Finally send $n\to \infty$, by \reff{BCD}  and applying the dominated convergence theorem under $\dbP$ and $\dbP^*$ we have
\beaa
\dbE^\dbP[ \wh Y_t \1_E] - \dbE^{\dbP^*}[\wh Y_{\wh\t^*-}1_E] &\le& \e.
\eeaa
That is,  $\dbP_\e := \dbP^*$ satisfies the requirement in the case $\t = t <\wh\t^*$ on $E$.
\\
{\it Step 2.} We now prove Claim \reff{ch2-claim}. 
Indeed,  for any $m\le n$ and any $\o\in E^{m,\d}_i$, by Lemma \ref{lem-YUComega} we have
\bea
\label{ch2-est1}
&&\wh Y_{t}(\o) - \dbE^{\dbP^{m,\d}_i}\Big[\wh Y^{t,\o}_{\bar\t_n^{t,\o}}\Big] \nonumber\\
&=&\wh Y_{t}(\o) - \wh Y_{t}(\o^{m,i}) + \wh Y_{t}(\o^{m,i}) -\dbE^{\dbP^{m,\d}_i}\Big[\wh Y^{t,\o^{m,i}}_{\bar\t_n^{t,\o^{m,i}}}\Big]  + \dbE^{\dbP^{m,\d}_i}\Big[\wh Y^{t,\o^{m,i}}_{\bar\t_n^{t,\o^{m,i}}} - \wh Y^{t,\o}_{\bar\t_n^{t,\o}}\Big] \nonumber\\ 
&\le& C\rho_1(\d) + \e +  \dbE^{\dbP^{m,\d}_i}\Big[\big|\wh Y^{t,\o^{m,i}}_{\bar\t_n^{t,\o^{m,i}}} - \wh Y^{t,\o}_{\bar\t_n^{t,\o}}\big| \1_{(\O^*_n)^{t, \o^{m,i}} \cap (\O^*_n)^{t,\o}} \psi_n^{t, \o^{m,i}}  \psi_n^{t,\o}\Big]  \\
&&+ C\dbP^{m,\d}_i\Big[[(\O^*_n)^{t, \o^{m,i}}]^c \cup [(\O^*_n)^{t,\o}]^c\Big] + C \dbE^{\dbP^{m,\d}_i}\Big[1-  \psi_n^{t, \o^{m,i}}  + 1-  \psi_n^{t, \o}\Big].\nonumber
\eea
Note that
\bea
\label{ch2-est2}
\dbE^{\dbP^{m,\d}_i}\Big[1-  \psi_n^{t, \o^{m,i}}  + 1-  \psi_n^{t, \o}\Big] &\le&2 \dbE^{\dbP^{m,\d}_i}\Big[1-  \psi_n^{t, \o}\Big] + \rho_n(\d);\nonumber\\
\dbP^{m,\d}_i\Big[[(\O^*_n)^{t, \o^{m,i}}]^c \cup [(\O^*_n)^{t,\o}]^c\Big]  &\le& 2\dbP^{m,\d}_i\Big[ [(\O^*_n)^{t,\o}]^c\Big]  +\dbP^{m,\d}_i\Big[[(\O^*_n)^{t, \o^{m,i}}]^c \cap (\O^*_n)^{t,\o} \Big]\\
&\le&2\dbP^{m,\d}_i\Big[ [(\O^*_n)^{t,\o}]^c\Big]  +  \dbP^{m,\d}_i\Big[0<d\big(\o \otimes_t B^t,  (\O^*_n)^c\big) < \d\Big]\nonumber\\
&\le& 2\dbP^{m,\d}_i\Big[ [(\O^*_n)^{t,\o}]^c\Big]  + \dbE^{\dbP^{m,\d}_i}\Big[h_\d \big(d\big(\o \otimes_t B^t,  (\O^*_n)^c\big)\big)\Big].\nonumber
\eea
Moreover, on $(\O^*_n)^{t, \o^{m,i}} \cap (\O^*_n)^{t,\o} \cap \{\psi_n^{t, \o^{m,i}}>0\} \cap\{ \psi_n^{t,\o}>0\}$, by  Lemma \ref{lem-Omegan*} and \reff{psinthn} we have
\beaa
(\underline\th^*_n)^{t,\o^{m,i}} \le \bar\t_n^{t,\o^{m,i}} \le (\overline\th^*_n)^{t,\o^{m,i}} < \ch_{n+4}^{t,\o^{m,i}};\q (\underline\th^*_n)^{t,\o} \le \bar\t_n^{t,\o} \le (\overline\th^*_n)^{t,\o} < \ch_{n+4}^{t,\o}.
\eeaa
Then 
\beaa
&&\big|\wh Y^{t,\o^{m,i}}_{\bar\t_n^{t,\o^{m,i}}} - \wh Y^{t,\o}_{\bar\t_n^{t,\o}}\big|\; \le\; \big|\wh Y^{t,\o^{m,i}}_{(\underline\th^*_n)^{t,\o^{m,i}}} - \wh Y^{t,\o}_{(\underline\th^*_n)^{t,\o}} \big|\\
&&\q+ \sup_{(\underline\th^*_n)^{t,\o^{m,i}} \le s \le (\overline\th^*_n)^{t,\o^{m,i}}}|\wh Y^{t,\o^{m,i}}_s- \wh Y^{t,\o^{m,i}}_{(\underline\th^*_n)^{t,\o^{m,i}}} | + \sup_{(\underline\th^*_n)^{t,\o} \le s \le (\overline\th^*_n)^{t,\o}}|\wh Y^{t,\o}_s- \wh Y^{t,\o}_{(\underline\th^*_n)^{t,\o}} | \\
& =&\big| \wh Y^{t,\o^{m,i}}_{(\underline\th^*_n)^{t,\o^{m,i}}} - \wh Y^{t,\o}_{(\underline\th^*_n)^{t,\o}}\big|+ 2 \sup_{(\underline\th^*_n)^{t,\o} \le s \le (\overline\th^*_n)^{t,\o}}|\wh Y^{t,\o}_s- \wh Y^{t,\o}_{(\underline\th^*_n)^{t,\o}} |  \\
&&+ \sup_{(\underline\th^*_n)^{t,\o^{m,i}} \le s \le (\overline\th^*_n)^{t,\o^{m,i}}}|\wh Y^{t,\o^{m,i}}_s- \wh Y^{t,\o^{m,i}}_{(\underline\th^*_n)^{t,\o^{m,i}}} | - \sup_{(\underline\th^*_n)^{t,\o} \le s \le (\overline\th^*_n)^{t,\o}}|\wh Y^{t,\o}_s- \wh Y^{t,\o}_{(\underline\th^*_n)^{t,\o}} |.
\eeaa
Applying Lemma \ref{lem-Ycontt} we get
\beaa
\big| \wh Y^{t,\o^{m,i}}_{(\underline\th^*_n)^{t,\o^{m,i}}} - \wh Y^{t,\o}_{(\underline\th^*_n)^{t,\o}}\big| &\le& Cn \rho_2\Big(\dbf_\infty\big(((\underline\th^*_n)^{t,\o^{m,i}}, \o^{m,i}\otimes_t B^t), ((\underline\th^*_n)^{t,\o}, \o\otimes_t B^t)\big)\Big)\\
&\le&  Cn \rho_2\Big(\d + \rho_n(\d) + 2\sup_{(\underline\th^*_n)^{t,\o}-\rho_n(\d) \le s \le  (\underline\th^*_n)^{t,\o} + \rho_n(\d)}|B^t_s - B^t_{(\underline\th^*_n)^{t,\o}}|\Big)\\
&\le&  Cn \rho_2\Big(\d + \rho_n(\d) + 2 \eta_n(\d)\Big),
\eeaa
and, similarly,
\beaa
&&\sup_{(\underline\th^*_n)^{t,\o^{m,i}} \le s \le (\overline\th^*_n)^{t,\o^{m,i}}}|\wh Y^{t,\o^{m,i}}_s- \wh Y^{t,\o^{m,i}}_{(\underline\th^*_n)^{t,\o^{m,i}}} | - \sup_{(\underline\th^*_n)^{t,\o} \le s \le (\overline\th^*_n)^{t,\o}}|\wh Y^{t,\o}_s- \wh Y^{t,\o}_{(\underline\th^*_n)^{t,\o}} | \\
&\le& \sup_{(\underline\th^*_n)^{t,\o^{m,i}} \le s \le(\underline\th^*_n)^{t,\o^{m,i}} \vee (\underline\th^*_n)^{t,\o}}|\wh Y^{t,\o^{m,i}}_s- \wh Y^{t,\o^{m,i}}_{(\underline\th^*_n)^{t,\o^{m,i}}} | \\
&&+ \sup_{(\underline\th^*_n)^{t,\o^{m,i}} \vee (\underline\th^*_n)^{t,\o} \le s \le (\overline\th^*_n)^{t,\o^{m,i}} \wedge(\overline\th^*_n)^{t,\o} }|\wh Y^{t,\o^{m,i}}_s- \wh Y^{t,\o}_s |+ | \wh Y^{t,\o^{m,i}}_{(\underline\th^*_n)^{t,\o^{m,i}}} - \wh Y^{t,\o}_{(\underline\th^*_n)^{t,\o}} | \\
&&+ \sup_{(\overline\th^*_n)^{t,\o^{m,i}} \wedge(\overline\th^*_n)^{t,\o} \le s \le(\overline\th^*_n)^{t,\o^{m,i}} }|\wh Y^{t,\o^{m,i}}_s- \wh Y^{t,\o^{m,i}}_{(\overline\th^*_n)^{t,\o^{m,i}}} |+| \wh Y^{t,\o^{m,i}}_{(\overline\th^*_n)^{t,\o^{m,i}}} - \wh Y^{t,\o}_{(\overline\th^*_n)^{t,\o}} |\\
&\le& Cn\rho_2\Big(\d + \rho_n(\d) + 2 \eta_n(\d)\Big) + C\rho_1(\d)\le Cn\rho_2\Big(\d + \rho_n(\d) + 2 \eta_n(\d)\Big).
\eeaa
Then
\beaa
\big|\wh Y^{t,\o^{m,i}}_{\bar\t_n^{t,\o^{m,i}}} - \wh Y^{t,\o}_{\bar\t_n^{t,\o}}\big| &\le&  Cn\rho_2\Big(\d + \rho_n(\d) + 2 \eta_n(\d)\Big) + 2 \sup_{(\underline\th^*_n)^{t,\o} \le s \le (\overline\th^*_n)^{t,\o}}|\wh Y^{t,\o}_s- \wh Y^{t,\o}_{(\underline\th^*_n)^{t,\o}} |.
\eeaa 
Plug this and \reff{ch2-est2} into \reff{ch2-est1},  for $\o\in E^{m,\d}_i$ we obtain
\beaa
&&\wh Y_{t}(\o) - \dbE^{\dbP^{m,\d}_i}\Big[\wh Y^{t,\o}_{\bar\t_n^{t,\o}}\Big] \;\le\; Cn \dbE^{\dbP^{m,\d}_i}\Big[\rho_2\Big(\d + \rho_n(\d) + 2 \eta_n(\d)\Big) \Big]+ C\rho_n(\d)+\e \\
&&\qq+ 2 \dbE^{\dbP^{m,\d}_i}\Big[\sup_{(\underline\th^*_n)^{t,\o} \le s \le (\overline\th^*_n)^{t,\o}}|\wh Y^{t,\o}_s- \wh Y^{t,\o}_{(\underline\th^*_n)^{t,\o}} |  \psi_n^{t,\o}\Big] \\
&&\qq  + C\dbP^{m,\d}_i\Big[ [(\O^*_n)^{t,\o}]^c\Big]+C \dbE^{\dbP^{m,\d}_i}\Big[1-  \psi_n^{t, \o}\Big] +C \dbE^{\dbP^{m,\d}_i}\Big[h_\d \big(d\big(\o \otimes_t B^t,  (\O^*_n)^c\big)\big)\Big].
\eeaa
Then by \reff{Pnd} we have, for any $N\ge n$, 
\bea
\label{ch2-est3}
&&\dbE^\dbP[ \wh Y_t \1_E] - \dbE^{\dbP^{N,\d}}[\wh Y_{t \vee \bar\t_n}\1_E] = \dbE^{\dbP^{N,\d}}\Big[[\wh Y_t -\wh Y_{\bar\t_n}]\1_{E\cap \{t< \bar\t_n\}}\Big] \nonumber\\
&\le& Cn\dbE^{\dbP^{N,\d}}\Big[\rho_2\Big(\d + \rho_n(\d) + 2 \eta_n(\d)\Big)\Big] + C\rho_n(\d)+\e+ C\dbP^{N,\d}\Big[ [\O^*_n]^c\Big]+C \dbE^{\dbP^{N,\d}}\Big[1-  \psi_n\Big]\nonumber\\
&&+ 2 \dbE^{\dbP^{N,\d}}\Big[\sup_{\underline\th^*_n \le s \le \overline\th^*_n}|\wh Y_s- \wh Y_{\underline\th^*_n} |  \psi_n\1_E\Big]    +C   \dbE^{\dbP^{N,\d}}\Big[h_\d\big(d\big(\o,  (\O^*_n)^c\big) \big)\Big]\nonumber\\
&\le&  Cn\cE\Big[\rho_2\Big(\d + \rho_n(\d) + 2 \eta_n(\d)\Big)\Big] + C\rho_n(\d)+\e + C 2^{-n} + C\cC(\psi_n < 1)\nonumber\\
&&+ 2 \dbE^{\dbP^{N,\d}}\Big[\sup_{\underline\th^*_n \le s \le \overline\th^*_n}|\wh Y_s- \wh Y_{\underline\th^*_n} |  \psi_n\1_E\Big] +C   \cE\Big[h_\d\big(d\big(\o,  (\O^*_n)^c\big) \big)\Big].
\eea
Similarly we have
\beaa
\dbE^{\dbP^{N,\d}}\Big[[\wh Y_{t \vee \bar\t_n} -\wh Y_{t \vee \underline\th^*_n}\psi_n] \1_E\Big]  &\le&  C 2^{-n} + C\cC(\psi_n < 1) + \dbE^{\dbP^{N,\d}}\Big[[\wh Y_{t \vee \bar\t_n} -\wh Y_{t \vee \underline\th^*_n}] \1_{E\cap \O^*_n} \psi_n\Big] \\
&\le&C 2^{-n} + C\cC(\psi_n < 1) +2\dbE^{\dbP^{N,\d}}\Big[\sup_{\underline\th^*_n \le s \le \overline\th^*_n}|\wh Y_s- \wh Y_{\underline\th^*_n} |  \psi_n\1_E\Big]
\eeaa
This, together with \reff{ch2-est3}, implies \reff{ch2-claim}.
\\
{\it Step 3.}  Finally we prove the lemma for general stopping time $\t$. We follow the arguments in Lemma \ref{lem-ch1}. Let $\t^n$ be a sequence of stopping times such that $\t^n \downarrow \t$ and each $\t^n$ takes only finitely many values. By applying the dominated convergence Theorem under $\dbP$, we may 
%\beaa
%\lim_{n\to\infty} \dbE^\dbP\Big[|\wh Y_{\t^n \wedge \wh\t^*} -\wh Y_\t|\1_E\Big]= 0. % \lim_{n\to\infty} \dbP\Big[E\cap \{\t^n\ge \wh\t^* > \t\}\Big]   =0 .
%\eeaa
% Let $\d>0$ be a constant which will be decided later. 
fix $n$ such that
\beaa
\dbE^\dbP\Big[|\wh Y_{\t^n \wedge \wh\t^*} -\wh Y_\t|\1_E\Big] \le {\e\over 2}. %,\q %\dbP\Big[E\cap \{\t < \wh\t^*\}\cap \tilde E^c\Big]= 
%\dbP\Big[E\cap \{\t^n\ge \wh\t^* > \t\}\Big] \le \d. % &\mbox{where}& \tilde E:= E\cap \{\t^n< \wh\t^* \} \in \cF_{\t^n}.
\eeaa
Assume $\t^n$ takes values $\{t_i, i=1,\cds, m\}$, and for each $i$,  denote $E_i := E\cap \{\t^n = t_i < \wh\t^*\}\in \cF_{t_i}$. Then $\{E_i, 1\le i\le m\}$ form a partition of $\tilde E := E\cap \{\t^n < \wh\t^*\}$. For each $i$, by Step 1 there exists $\dbP^i \in \cP(\dbP, t_i, E_i)$ such that
\beaa
\dbE^{\dbP}[\wh Y_{t_i} \1_{E_i}\Big] \le  \dbE^{\dbP^i}[\wh Y_{\wh\t^*-} \1_{E_i}\Big] + {\e \over 2m}.
\eeaa
Now define $\dbP_\e := \sum_{i=1}^m \dbP^i\1_{E_i} + \dbP \1_{\tilde E^c} \in \cP(\dbP, \t^n, \tilde E)\subset \cP(\dbP, \t, E)$.  Recall that $\tilde E\in \cF_{\t^n}$ and note that $\wh Y_{\wh\t^*} \le \wh Y_{\wh\t^*-}$, thanks to the supermartingale property of $\wh Y$. Then
\beaa
\dbE^\dbP\Big[ \wh Y_{\t} \1_E\Big]- \dbE^{\dbP_\e}\Big[\wh Y_{\wh\t^*-}\1_E\Big] &\le& {\e\over 2} +\dbE^\dbP\Big[ \wh Y_{\t^n\wedge \wh\t^*} \1_E\Big]- \dbE^{\dbP_\e}\Big[\wh Y_{\wh\t^*-}\1_E\Big] \\
&\le&{\e\over 2} +  \dbE^\dbP\Big[ \wh Y_{\t^n} \1_{\tilde E}\Big]- \dbE^{\dbP_\e}\Big[\wh Y_{\wh\t^*-}\1_{\tilde E}\Big] \\
&=&{\e\over 2} + \sum_{i=1}^m\Big(  \dbE^\dbP\Big[ \wh Y_{t_i} \1_{E_i}\Big]- \dbE^{\dbP_\e}\Big[\wh Y_{\wh\t^*-}\1_{E_i}\Big]\Big) \\
&\le&{\e\over 2}  + \sum_{i=1}^m{\e\over 2m}  = \e.
\eeaa 
The proof is complete now.
\qed

We need one more lemma.

\begin{lem}
\label{lem-ch3}
Let $\dbP\in \cP$, $\t\in \cT$, and $E\in \cF_\t$  such that $\t \le \ch$ on $E$. For any  $\e>0$, there exists $\dbP_\e\in\cP(\dbP,\t, E)$ such that 
\beaa
 \ch  &\le& \t +  {1\over L}d(\o_\t, O^c) + 3\e + \sup_{\t\le t\le\t+ \e} |\o_t-\o_\t|, ~~\dbP_\e\mbox{-a.s. on}~E 
\eeaa
\end{lem}
\proof   First, there exists $\tilde\t \in \cT$ such that $\t \le \tilde\t \le \t + \e$ and $\tilde \t$ takes only finitely many values $0\le t_1<\cds<t_n = t_0$. Denote $E_i :=E\cap  \{\tilde \t = t_i <\ch\}\in \cF_{t_i}$. Then $\{E_i, 1\le i\le n\}$ is a partition of $E\cap \{\tilde\t <\ch\}$ and
\bea
\label{chtau1}
\ch  \le \tilde\t \le \t+ \e &\mbox{on}& E\cap \{\tilde\t \ge\ch\}.
\eea

  For any $i$, there exists a partition $(E^i_j)_{j\ge 1}$ of $E_i$ such that $|\o_{t_i}-\o'_{t_i}|\le L\e$ for any $\o, \o'\in E^i_j$. For each $(i, j)$, fix an $\o^{ij} \in E^i_j$ and a unit vector $\a^{ij}$ pointing to the direction from  $\o^{ij}_{t_i}$ to $O^c$. Now for any $\o\in E^i_j$, define $\dbP^{i,j,\o}\in \cP_{t_i}$ as follows: 
\beaa
\b = {\bf 0},  \q \a_t = {1\over \e}[\o^{ij}_{t_i} - \o_{t_i}] \1_{[t_i, t_i + \e)}(t) + L \a^{ij}\1_{ [t_i + \e, T]}(t).
\eeaa
We see that
\beaa
\ch^{t_i,\o} = \Big[t_i + \e + {1\over L}d(\o^{ij}_{t_i}, O^c)\Big]\wedge t_0, ~\dbP^{i,j,\o}\mbox{-a.s. on}~E^i_j.
\eeaa
Similar to the proof of \reff{DPP1claim}, there exists $\dbP_\e\in\cP(\dbP,\tilde \t, E)\subset \cP(\dbP, \t, E)$ such that  the r.c.p.d. $\dbP^{t_i,\o}_\e = \dbP^{i,j,\o}$ for $\dbP$-a.e. $\o\in E^i_j$. Then
\beaa
%\label{chtau2}
\ch &\le& \t + 2 \e + {1\over L}[d(\o_{t_i}, O^c) + L\e] \le  \t + 3\e + {1\over L}\Big[d(\o_\t, O^c) + |\o_\t-\o_{t_i}|\Big]\nonumber\\
&\le& \t + 3\e + {1\over L}\Big[d(\o_\t, O^c) +\sup_{\t\le t\le \t+\e} |\o_t-\o_{\t}|\Big] , \q \dbP_\e\mbox{-a.s. on }~E^i_j. 
\eeaa
This, together with \reff{chtau1}, proves the lemma.
\qed

We are now ready to complete the

\ms

\no{\bf Proof of Proposition \ref{prop-chcont}.} The inequality $\cE[\wh Y_{\wh\tau^*}]\le\cE[\wh Y_{\wh\tau^*-}]$ is a direct consequence of the $\cE-$supermartingale property of $\wh Y$ established in Theorem \ref{thm-optimal-DP}. As for the reverse inequality, 
since $\wh Y$ is  continuous on $[0,\ch)$ and $\ch_n \uparrow \ch$ with $\ch_n < \ch$, 
it suffices to show that, for any $\dbP\in \cP$ and any  $\e>0$
\bea
\label{chcont}
I_n := \dbE^\dbP[\wh Y_{\wh\t^* \wedge \ch_n}]  - 
       \cE [ \wh Y_{\wh\t^*}] \le 5\e
 &\mbox{for sufficiently large}&
 n.
\eea
Let $\d > 0$, $n> {1\over L\d}$.  Set $t_n := t_0-{1\over n}$, $\t^0 :=\wh \t^*\wedge \ch_n$, and $\dbP^0:=\dbP$. We proceed in two steps.

\ms

\no{\it Step 1.} Apply Lemma \ref{lem-ch1} with $\dbP^0, \t^0, {\wh\tau^*}$,  and $\O$, there exist $\dbP^{1,1}\in \cP(\dbP^0, \t^0, \O)$ and a stopping time $\tilde\t^1$ taking values in $[\t^0, {\wh\tau^*}]$,  such that
  \beaa
\dbE^{\dbP^0}[ \wh Y_{\t^0}]   \le \dbE^{\dbP^{1,1}}\Big[\wh X_{\tilde\t^1}  \1_{\{\tilde\t^1<{\wh\tau^*}\}} +\wh Y_{\wh\tau^*}\1_{\{\tilde\t^1={\wh\tau^*}\}}\Big] + \e.
\eeaa
Denote $E_1 := \{ \tilde\t^1 < t_n\}\in \cF_{\tilde\t^1}$. By \reff{USC2} and following the same argument as for the estimate in \reff{rho1est}, we have:  $\dbP^{1,1}\mbox{-a.s. on}~  E^c_1 \cap \{\tilde\t^1 <{\wh\tau^*}\}$,
\beaa
\wh X_{\tilde\t^1}  &\le& \wh X_{\tilde\t^1} - \dbE^{\dbP^{1,1}}_{\tilde\t^1}[\wh X_{\wh\tau^*}]  +  \dbE^{\dbP^{1,1}}_{\tilde\t^1}[\wh Y_{\wh\tau^*}] \\
&\le& \dbE^{\dbP^{1,1}}_{\tilde\t^1}\Big[\rho_0\big({1\over n} + \|B^{\tilde \t^1}\|_{\tilde\t^1+{1\over n}}\big)\Big]  + \dbE^{\dbP^{1,1}}_{\tilde\t^1}[\wh Y_{\wh\tau^*}] 
\le C\bar\rho_0({n^{-1}})+ \dbE^{\dbP^{1,1}}_{\tilde\t^1}[\wh Y_{\wh\tau^*}].
\eeaa
Then, denoting $E_2 :=  E_1 \cap \{\tilde\t^1 < {\wh\tau^*}\}\in \cF_{\tilde\t^1}$, we get:
  \bea
   \label{chest1}
 \dbE^{\dbP^{0}}\Big[\wh Y_{\t^0} \Big]  &\le& \dbE^{\dbP^{1,1}}\Big[\wh X_{\tilde\t^1} \1_{E_2} + \wh X_{\tilde\t^1}  \1_{E^c_1\cap \{\tilde\t^1<{\wh\tau^*}\}} +\wh Y_{\wh\tau^*}\1_{\{\tilde\t^1={\wh\tau^*}\}}\Big] + \e\nonumber\\
 &\le&  \dbE^{\dbP^{1,1}}\Big[\wh X_{\tilde \t^1}\1_{E_2} +  \wh Y_{{\wh\tau^*}} \1_{E_2^c}
                   \Big]  + C\bar\rho_0(n^{-1}) \dbP^0[E_1^c]  + \e.
                   \eea     

Next,  set $\tilde\d := [\d^2\bar\rho_0(3\d)]\wedge {\d\over 3}$. Apply Lemma \ref{lem-ch3} on $\dbP^{1,1}$, $\tilde\t^1$, $E_2$, and $\tilde\d$,  there exists $\dbP^{1,2} \in \cP(\dbP^{1,1}, \tilde\t^1, E_2)$ such that 
\beaa
%\dbP^{1,1}(E_3\backslash E_3^\e) \le \e &\mbox{and}& {\wh\tau^*}\le \ch \le \tilde \t^1 + 2\d,\q \dbP^{1,2}\mbox{-a.s. on}~E_3^\e\cap \{ d(\o_{\tilde\t^1}, O^c) \le L\d\}.
 \ch  &\le& \tilde\t^1 +  {1\over L}d(\o_{\tilde\t^1}, O^c) + \d+ \|\o^{\tilde\t^1}_t\|_{\tilde\t^1+\tilde\d}, ~~\dbP^{1,2}\mbox{-a.s. on}~E_2.
\eeaa
Since $\tilde\t^1 \le \wh\t^*\le \ch$, we have
\beaa
\wh\t^* - \tilde\t^1 \le 3\d, ~~ \dbP^{1,2}\mbox{-a.s. on}~E_2\cap \{ d(\o_{\tilde\t^1}, O^c) \le L\d\} \cap \{\|\o^{\tilde\t^1}\|_{\tilde\t^1 + \tilde\d} \le \d\}.
\eeaa
Then,  by \reff{USC2} and \reff{rho1est} again  we have: $ \dbP^{1,2}\mbox{-a.s. on}~E_2\cap \{ d(\o_{\tilde\t^1}, O^c) \le L\d\}\in \cF_{\tilde\t^1}$,
\beaa
\wh X_{\tilde\t^1} &\le& \dbE^{\dbP^{1,2}}_{\tilde\t^1} [\wh X_{{\wh\tau^*}}] + \dbE^{\dbP^{1,2}}_{\tilde\t^1}\Big[\rho_0\Big(\dbf_\infty\big((\tilde\t^1, B), (\wh\t^*, B) \big)\Big]\\
&=& \dbE^{\dbP^{1,2}}_{\tilde\t^1} [\wh X_{{\wh\tau^*}}] + \dbE^{\dbP^{1,2}}_{\tilde\t^1}\Big[\rho_0\Big(\dbf_\infty\big((\tilde\t^1, B), (\wh\t^*, B) \big)\big[\1_{\{\|B^{\tilde\t^1}\|_{\tilde\t_1+\tilde\d} \le \d\} } + \1_{\{\|B^{\tilde\t^1}\|_{\tilde\t_1+\tilde\d}> \d\} }\big] \Big]\\
&\le& \dbE^{\dbP^{1,2}}_{\tilde\t^1} [\wh X_{{\wh\tau^*}}] + \dbE^{\dbP^{1,2}}_{\tilde\t^1}\Big[\rho_0\Big(3\d + \|B^{\tilde\t^1}\|_{\tilde\t^1+3\d}\Big)\Big] + C\d^{-2}  \dbE^{\dbP^{1,2}}_{\tilde\t^1}[ \|B^{\tilde\t^1}\|_{\tilde\t_1+\tilde\d}^2]\\
&\le& \dbE^{\dbP^{1,2}}_{\tilde\t^1} [\wh X_{{\wh\tau^*}}] + C\bar\rho_0(3\d) + {C\tilde\d\over \d^2}\le  \dbE^{\dbP^{1,2}}_{\tilde\t^1} [\wh X_{{\wh\tau^*}}] + C\bar\rho_0(3\d).
\eeaa
Note that $n^{-1} \le L\d \le 3\d$. Thus, denoting $E_3 :=  E_2 \cap \{ d(\o_{\tilde\t^1}, O^c) > L\d\}\in \cF_{\tilde\t^1}$, \reff{chest1} leads to:
  \bea
  \label{chest2}
\dbE^{\dbP^{0}}\Big[\wh Y_{\t^0} \Big]  
                    &\le& \dbE^{\dbP^{1,2}}\Big[\wh X_{\tilde \t^1}\1_{E_3}+  \wh Y_{{\wh\tau^*}} \1_{E_3^c}
                   \Big]   +C\bar\rho_0(3\d)\dbP^{1,2}(E_3^c)+ \e.
                   \eea
                   
Moreover, apply  Lemma \ref{lem-ch2} with $\dbP^{1,2}, \tilde\t^1$, $E_3$, and $\e$,  there exists $\dbP^{1,3}\in \cP(\dbP^{1,2}, \tilde\t^1, E_3)$ such that
\beaa
\dbE^{\dbP^{1,2}}\Big[\wh X_{\tilde\t^1}\1_{E_3}\Big] \le \dbE^{\dbP^{1,2}}\Big[ \wh Y_{\tilde\t^1}\1_{E_3}\Big] \le \dbE^{\dbP^{1,3}}\Big[\wh Y_{\wh\t^*-}\1_{E_3}\Big] + \e.
\eeaa
Define $\t^1 := \inf\{t\ge \tilde\t^1: d(\o_t, O^c) \le {1\over n}\} \wedge {\wh\tau^*}$. Note that $\t^1 < \ch$ on $E_3$ and $\wh Y$ is a $\dbP^{1,3}$-supermartingale. Then
\beaa
\dbE^{\dbP^{1,3}}\Big[\wh Y_{\wh\t^*-}\1_{E_3}\Big] &\le& \dbE^{\dbP^{1,3}}\Big[\wh Y_{\t^1}\1_{E_3}\Big].
\eeaa  
Thus
\beaa
\dbE^{\dbP^{1,2}}\Big[\wh X_{\tilde\t^1}\1_{E_3}\Big] &\le& \dbE^{\dbP^{1,3}}\Big[\wh Y_{\t^1}\1_{E_3}\Big] + \e.
\eeaa
Plug this into \reff{chest2}, we obtain
\beaa
\dbE^{\dbP^{0}}\Big[\wh Y_{\t^0} \Big]  
                    &\le& \dbE^{\dbP^{1,3}}\Big[\wh Y_{\t^1}\1_{E_3}+  \wh Y_{{\wh\tau^*}} \1_{E_3^c}
                   \Big]   +C\bar\rho_0(3\d)\dbP^{1,3}(E_3^c)+ 2\e.
                   \eeaa
We now denote $\dbP^1 := \dbP^{1,3}\in \cP(\dbP^0, \t^0, \O)$, and 
\bea
\label{D1}
D_1 := E_3 \cap \{\t^1 < {\wh\tau^*}\}=\{\tilde\t^1<t_n \wedge \tilde\t^*\} \cap \{d(\o_{\tilde\t^1}, O^c) > L\d\} \cap  \{\t^1 < {\wh\tau^*}\} \in\cF_{\t^1} 
\eea
 Then
 \bea
 \label{chest3}
\dbE^{\dbP^{0}}\Big[\wh Y_{\t^0} \Big]  
                    &\le& \dbE^{\dbP^{1}}\Big[\wh Y_{\t^1}\1_{D_1}+  \wh Y_{{\wh\tau^*}} \1_{D_1^c}
                   \Big]   +C\bar\rho_0(3\d) \dbP^{1}(D_1^c)+ 2\e.
                   \eea
                   
                   \ms
\no{\it Step 3:} Iterating the arguments of Step 1, we may define $(\tilde\t^m, \t^m, \dbP^m, D_m)_{m\ge 1}$ such that:
 \beaa
 &&  \dbP^{m+1}\in \cP(\dbP^{m}, \t^{m}, D_m),\q \t^{m} \le \tilde\t^{m+1} \le {\wh\tau^*};\\
 && \t^{m+1} := \inf\Big\{ t \ge \tilde\t^{m+1} : d(\o_t, O^c) \le {1\over n}\Big\} \wedge {\wh\tau^*}\\
 &&  D_{m+1} := D_{m}\cap \{\tilde\t^{m+1} < t_n\wedge {\wh\tau^*}\} \cap   \{d(\o_{\tilde\t^{m+1}}, O^c) > L\d\} \cap \{\t^{m+1} <{\wh\tau^*}\};
  \eeaa 
 and
 \beaa
 \dbE^{\dbP^m}\Big[\wh Y_{\t^{m}}\1_{D_{m}} \Big] 
 &\le& 
 \dbE^{\dbP^{m+1}}
 \Big[ \wh Y_{\t^{m+1}}\1_{D_{m+1}} + \wh Y_{{\wh\tau^*}}\1_{D_m\cap D^c_{m+1}} 
       \Big] 
 \\
 &&+C\bar\rho_0(3\d)\dbP^{m+1}(D_m\cap D^c_{m+1})+ 2^{1-m}\e.
 \eeaa
By induction, for any $m\ge 1$ we have
 \bea
 \label{chest4}
\dbE^{\dbP^{0}}\Big[\wh Y_{\t^0} \Big]  
                    &\le& \dbE^{\dbP^{m}}\Big[\wh Y_{\t^m}\1_{D_m}+  \wh Y_{{\wh\tau^*}} \1_{D_m^c}
                   \Big]   +C\bar\rho_0(3\d) \dbP^{m}(D_m^c)+ 4\e\nonumber\\
                   &\le&  \dbE^{\dbP^{m}}[\wh Y_{{\wh\tau^*}}] + 2C_0 \dbP^m[D_m]  +C\bar\rho_0(3\d) + 4\e.
                   \eea
Note that
 \beaa
 \dbP^m[D_m] &\le&  \dbP^m\Big[\cap_{i=1}^m \{|B_{\tilde\t^i} - B_{\t^{i-1}}| \ge L\d - {1\over n}\} \cap  \{|B_{\t^i} - B_{\tilde\t^{i}}| \ge L\d - {1\over n}\} \Big] \\
 &\le&  \dbP^m\Big[\sum_{i=1}^m [|B_{\tilde\t^i} - B_{\t^{i-1}}|^2 +  |B_{\t^i} - B_{\tilde\t^{i}}|^2] \ge 2m(L\d - {1\over n})^2 \Big] \\
 &\le& {1\over 2m(L\d - {1\over n})^2} \dbE^{\dbP^m}\Big[\sum_{i=1}^m [|B_{\tilde\t^i} - B_{\t^{i-1}}|^2 +  |B_{\t^i} - B_{\tilde\t^{i}}|^2\Big] \le {C\over 2m(L\d - {1\over n})^2}.
 \eeaa
Then, \reff{chest4} leads to
  \beaa
I_n \le  {C\over 2m(L\d - {1\over n})^2} + C\bar\rho_0(3\d)+ 4\e.
 \eeaa
which implies, by sending $m\to \infty$ that
 \beaa
I_n \le  C\bar\rho_0(3\d) + 4\e.
 \eeaa 
Hence, by choosing $\d$ small enough such that $\bar\rho_0(3\d) \le \e$, we see that \reff{chcont} holds true for $n> {1\over L\d}$.
 \qed

\end{document}